\documentclass{article}
\usepackage{amsmath}
\usepackage{amsfonts}
\usepackage{amssymb}
\usepackage{amsthm}
\usepackage{amscd}
\usepackage{makeidx}
\makeindex
\newtheorem{theorem1}{Theorem}[section]
\newtheorem{remark}{Remark}[section]

\newtheorem{proposition1}{Proposition}[section]

\newtheorem{definition1}{Definition}[section]

\newtheorem{lemma1}{Lemma}[section]

\newtheorem{cor1}{Corollary}[section]

\newcommand{\supp}{\text{supp }}
\newcommand{\e}{\varepsilon}
\newcommand{\M}{{\mathcal M}}
\newcommand{\R}{{\mathcal R}}
\begin{document}

\title{Operator synthesis. I.\\
Synthetic sets, bilattices and tensor algebras.}

\author{Victor Shulman \vspace{0.1cm}\\
{\footnotesize \sl Department of Mathematics, Vologda State Liceum of
Mathematical and Natural Sciences,
}\\
{\footnotesize \sl
Vologda, 160000, Russia} \vspace{0.3cm}\\
Lyudmila  Turowska \footnote{Corresponding author.
E-mail address: turowska\symbol{64}math.chalmers.se. Fax: +46 31 161973} \vspace{0.1cm}\\
{\footnotesize \sl Department of Mathematics, Chalmers University of 
Technology, }\\ 
{\footnotesize \sl SE-412 96 G\"oteborg, Sweden} }

\date{}

\maketitle{}
\begin{abstract}
The interplay between the invariant subspace theory and spectral synthesis for 
locally compact abelian group discovered by Arveson \cite{arv} is extended to 
include other topics as harmonic analysis for Varopoulos algebras and 
approximation by projection-valued measures. 
We propose a ''coordinate'' approach which nevertheless does not use the technique of
pseudo-integral operators, as well as a coordinate free one which allows to 
extend to non-separable spaces some important results and constructions of 
\cite{arv}. We solve some problems posed in \cite{arv}.
\end{abstract}

\vspace{1cm}


\section{Introduction}
   The classical notion of spectral synthesis is related to the Galois
correspondence between ideals $J$ of a commutative regular Banach 
algebra ${\mathcal A}$
and closed subsets $E$ of its character space $X({\mathcal A})$:
$\ker J = \{t \in X({\mathcal A}): t(a) = 0, \text{ for any }a \in J\}$,
$\text{hull }E = \{a \in {\mathcal A}: t(a) = 0, \text{ for any }t \in E\}$.
Namely, a set $E$ is called synthetic (or a set of spectral synthesis)
if $\ker J = E$ implies $J = \text{hull }E$. Note, that the converse 
implication holds for any closed $E \subseteq X({\mathcal A})$.

In the invariant subspace theory
the central object is a Galois correspondence between 
operator algebras ${\mathcal M}$ and strongly closed subspace 
lattices $\mathcal L$:
$\text{lat }{\mathcal M} = \{L: TL \subseteq L, \text{for any }T \in 
{\mathcal M}\}$,
$\text{alg }{\mathcal L} = \{T: TL \subseteq L, \text{for any }L \in 
{\mathcal L}\}$.
A lattice ${\mathcal L}$ can be called {\it operator synthetic} if 
$\text{lat }{\mathcal M} = {\mathcal L}$ implies ${\mathcal M} = 
\text{alg }{\mathcal L}$.

W.Arveson \cite{arv} proved that if one restricts the map lat to the 
variety of algebras, containing a fixed maximal abelian selfadjoint
algebra (masa), then the above formal analogy becomes very rich and 
fruitful. In particular, answering a question of H.Radjavi and
P.Rosenthal, he proved the failure of operator synthesis
in the class of $\sigma$-weakly closed algebras, containing a masa 
(Arveson algebras, in terminology of \cite{EKS}), by using the famous 
L.Schwartz's example of a non-synthetic set for the group algebra 
$L^1({\mathbb R}^3)$. Note, that among other brilliant results, 
\cite{arv} contains the implication 
${\mathcal M} = \text{alg }{\mathcal L}\Rightarrow{\mathcal L} = 
\text{lat }{\mathcal M}$, for an Arveson algebra ${\mathcal M}$ 
(in full analogy with the classical situation).

The results in  \cite{arv} indicate, in fact, that the problematic of 
the operator synthesis obtains a more natural setting if instead of 
algebras and lattices one considers bimodules over masas and their
bilattices (see the definitions below). We choose this point of view
aiming at the investigation of various faces of operator synthesis,
that reflect its connections with measure theory, approximation theory,
linear operator equations and spectral theory of multiplication 
operators, synthesis in modules, Haagerup tensor products and Varopoulos
tensor algebras.

   Let us list some results, proved in this first part of 
our work. We show the equivalence of several different definitions
of operator synthesis. Answering a question of W.Arveson 
we prove the existence of a minimal pre-reflexive algebra (bimodule) with a
given invariant subspace lattice (bilattice), without the assumption 
of separability of the underlying Hilbert space. On the other hand, for 
separable case we propose a coordinate approach which does not need a choice 
of a topology, replacing it by the pseudo-topology, naturally related to 
the measure spaces. This allows to consider simultaneously the synthesis 
for a more wide class of subsets and to avoid the use of pseudo-integral 
operators and the complicated theory of integral decompositions of 
measures (see \cite{arv} and \cite{davidson}). This approach admits also 
the use of measurable sections which leads to an "inverse image theorem" 
(Theorem~\ref{preimage}) for operator synthesis, implying in particular 
Arveson's theorem on synthesis for finite width lattices. We answer 
(in the negative) a question posed by Arveson \cite{arv}[Problem, p.487] 
on synthesizability of the lattice generated by a synthetic lattice and a 
lattice of finite width (Theorem~\ref{synfw}).We prove that a closed 
subset in a product of two compact sets is a set of spectral 
synthesis for the Varopoulos algebra if it is operator synthetic for 
any choice of measures (Theorem~\ref{v})(Proposition~\ref{failure} shows
that the converse implication fails). This, together with the above mentioned
inverse image theorem, gives some sufficient conditions for spectral 
synthesis, implying, for example, the well known Drury's theorem on 
non-triangular sets (Corollary~\ref{finw}).

   In the second part of the work we are going to consider the
individual operator synthesis and its connections with linear operator
equations.

 We are indebted to S.Drury, A.Katavolos, S.Kaijser,
B.Magajna, I.Todorov, N.Varopou\-los, 
for helpful discussions and valuable information. We thank the referee for
several important suggestion.
The work was partially written when the first author was visiting 
Chalmers University of Technology in G\"oteborg, Sweden.
The research  was partially supported by a grant from the Swedish Royal 
Academy of Sciences as a part of the program of cooperation with former
Soviet Union.

\section{Synthetic sets (measure-theoretic approach)}\label{section2}
Let $(X,\mu)$, $(Y,\nu)$ denote $\sigma$-finite separable spaces with 
standard measures.We use standard measure-theoretic terminology. A subset 
of the Cartesian product
$X\times Y$ is said to be a measurable rectangle if it has the form 
$A\times B$ with measurable $A\subseteq X$, $B\subseteq Y$. A set 
$E\subseteq X\times Y$ is called {\it marginally null set} if
$E\subseteq(X_1\times Y)\cup(X\times Y_1)$, where $\mu(X_1)=\nu(Y_1)=0$.
If subsets $\alpha$, $\beta$ of $X\times Y$ are marginally equivalent 
(i.e. their symmetric difference is marginally null) we write 
$\alpha\cong\beta$.
Following \cite{EKS} we define $w$-topology on $X\times Y$ such that the
$w$-open (pseudo-open) sets are, modulo marginally null sets, countable
union of measurable rectangles. The complements of $w$-open sets are called 
$w$-closed (pseudo-closed). The complement to a set $A$ will be denoted by
$A^c$. 

Let $\Gamma(X,Y)=L_2(X,\mu)\hat\otimes L_2(Y,\nu)$ be the projective tensor
product, i.e. the space of all functions $F:X\times Y\to{\mathbb C}$ which
admit a representation
\begin{equation}\label{eq1}
F(x,y)=\sum_{n=1}^{\infty}f_n(x)g_n(y)
\end{equation}
where $f_n\in L_2(X,\mu)$, $g_n\in L_2(Y,\nu)$   and 
$\sum_{n=1}^{\infty}||f_n||_{L_2}\cdot||g_n||_{L_2}<\infty$. Such a function
$F$ is defined  marginally almost everywhere (m.a.e.) in that,  if 
$f_n$, $g_n$ are changed on null sets then $F$ will change on 
a marginally null set. Then 
$L_2(X,\mu)\hat\otimes L_2(Y,\nu)$-norm of such a function $F$ is
$$||F||_{\Gamma}=\inf\sum_{n=1}^{\infty}||f_n||_{L_2}\cdot
||g_n||_{L_2},$$ 
where the infinum is taken over all  sequences $f_n$, $g_n$ for which
(\ref{eq1}) holds m.a.e. In what follows we 
identify two functions in  $\Gamma(X,Y)$ which coincides m.a.e.

By \cite{EKS}[Theorem~6.5], any function $F\in \Gamma(X,Y)$ is 
pseudo-continuous (continuous with respect to the $\omega$-topology defined 
above). We say that $F\in \Gamma(X,Y)$ vanishes
on a  set $K\subseteq X\times Y$ if $F\chi_K=0$ (m.a.e), where
$\chi_K$ is the characteristic function of  $K$.
For arbitrary $K\subseteq X\times Y$ denote by $\Phi(K)$ the set of 
all functions
$F\in\Gamma(X,Y)$ vanishing on $K$. Clearly $\Phi(K)$ is a subspace
of $\Gamma(X,Y)$.
\begin{lemma1}\label{lemma1}
Any convergent in norm sequence $\{F_n\}\in\Gamma(X,Y)$ has a subsequence
which converges marginally almost everywhere.
\end{lemma1}
\begin{proof}
We may assume that $\{F_n\}$ converges to zero in norm. Then there exist
 functions $f_k^{(n)}\in L_2(X,\mu)$, $g_k^{(n)}\in L_2(Y,\nu)$ such that
$$F_n(x,y)=\sum_{k=1}^{\infty}f_k^{(n)}(x)g_k^{(n)}(y),\quad
\sum_{k=1}^{\infty}||f_k^{(n)}||^2_{L_2}\to 0 \text{ and }
\sum_{k=1}^{\infty}||g_k^{(n)}||^2_{L_2}\to 0.$$ By the Riesz theorem
applied to the functions $f^{(n)}(x)=\sum_{k=1}^{\infty}|f_k^{(n)}(x)|^2$
and $g^{(n)}(y)=\sum_{k=1}^{\infty}|g_k^{(n)}(y)|^2$ there exists a
subsequence $\{F_{n_j}\}$ such that
$f^{(n_j)}(x)$ and $g^{(n_j)}(y)$ converge to zero almost everywhere.
Therefore,  there exist $M\subset X$, $N\subset Y$, $\mu(M)=0$, $\nu(N)=0$,
such that
$f^{(n_j)}(x)\to 0$ and $g^{(n_j)}(y)\to 0$ for any $x\in X\setminus M$,
$y\in Y\setminus N$,  and since 
$|F_{n_j}(x,y)|\leq  f^{(n_j)}(x)g^{(n_j)}(y)$, this implies 
$F_{n_j}(x,y)\to 0$ for any $(x,y)\in (X\setminus M)\times (Y\setminus N)$.
\end{proof}
\begin{proposition1}\label{cl} $\Phi(K)$ is closed.
\end{proposition1}
\begin{proof}
Let $F\in\overline{\Phi(K)}$. By Lemma~\ref{lemma1} there exists a
sequence $F_n\in\Phi(K)$ which converges to $F$ marginally almost everywhere.
Removing a countable union of marginally null sets we can assume
that all $F_n$ vanish on the rest of the set $K$ and therefore $F\chi_K=0$ 
m.a.e. 
\end{proof}

If $F\in\Gamma(X,Y)$ vanishes on $K$ then by pseudo-continuity it vanishes
on the pseudo-closure of
$K$ so that  without loss of generality we can restrict ourselves to 
pseudo-closed sets $K$.

Given arbitrary subset ${\mathcal F}\subseteq \Gamma(X,Y)$, we define the
null set of ${\mathcal F}$, $\text{null }{\mathcal F}$, to be the largest, 
up to marginally null sets, pseudo-closed set such that each function 
$F\in{\mathcal F}$
vanishes on it. To see the existence of such a set take a countable dense
subset ${\mathcal A}\subseteq{\mathcal F}$ and consider 
$K=\cap_{F\in{\mathcal A}}F^{-1}(0)$. Clearly, $K$ is pseudo-closed,
${\mathcal A}\subseteq \Phi(K)$ and, by Proposition~\ref{cl},
${\mathcal F}=\overline{\mathcal A}\subseteq\Phi(K)$. The maximality of
$K$ is obvious.

Let $\Phi_0(K)$ be the closure in $\Gamma(X,Y)$ of the set of all functions
which vanish on neighbourhoods of $K$ (pseudo-open sets containing $K$). 
$\Phi_0(K)$ is a closed subspace
of $\Phi(K)$.
\begin{proposition1}\label{null}
$\text{\rm null }\Phi_0(K)=K=\text{\rm null }\Phi(K)$.
\end{proposition1}    
\begin{proof}
We work modulo marginally null sets. Let $\alpha\subseteq X$, 
$\beta\subseteq Y$ be measurable sets 
such that $(\alpha\times\beta)\cap K=\emptyset$. Then the  function 
$\chi_{\alpha}(x)\chi_{\beta}(y)$ belongs to $\Phi_0(K)$ and therefore 
$\text{null }\Phi_0(K)\subseteq (\alpha\times\beta)^c$.
Since $K$ is pseudo-closed, $K=(\cup_{k=1}^{\infty}\alpha_k\times\beta_k)^c$ 
for some measurable $\alpha_k$, $\beta_k$ so that  
$(\alpha_k\times\beta_k)\cap K=\emptyset$ and thus 
$\text{null }\Phi_0(K)\subseteq K$.
We have also that $\text{null }\Phi_0(K)\supseteq \text{null }\Phi(K)\supseteq K$ 
which implies our result.
\end{proof}
Clearly, the subspaces $\Phi_0(K)$ and $\Phi(K)$ are invariant with respect 
to the multiplication by functions $f\in L_{\infty}(X,\mu)$ and 
$g\in L_{\infty}(Y,\nu)$ (we just write invariant).
\begin{theorem1}\label{th1}
If  $A\subseteq \Gamma(X,Y)$ is an invariant  closed subspace then 
\begin{equation}\label{incl}
\Phi_0(\text{\rm null }A)\subseteq A
\subseteq\Phi(\text{\rm null A}).
\end{equation}
\end{theorem1}
The second inclusion  is obvious.
The proof of the first one is postponed  till Section~\ref{sec4}. 
This theorem justifies the following definition. 
\begin{definition1} We say that a pseudo-closed set 
$K\subseteq X\times Y$ is
synthetic (or $\mu\times\nu$-synthetic) if 
$$\Phi_0(K)=\Phi(K).$$
\end{definition1}
We shall also refer to synthetic sets as sets of operator synthesis or sets
of $\mu\times\nu$-synthesis when the measures need to be specified.

We shall see that sets of operator synthesis can be defined in several 
different ways. The relation to operator theory is based on the fact
that elements of $\Gamma(X,Y)$ are the kernels of the nuclear 
(trace class) operators from
$H_2=L_2(Y,\nu)$ to $H_1=L_2(X,\mu)$ and  the space  
${\mathfrak S}^1(H_2,H_1)$ of all such operators is isometrically 
isomorphic to $\Gamma(X,Y)$ (see \cite{arv}).
The space of bounded operators, $B(H_1,H_2)$,  from 
$H_1$ to $H_2$ is dual to ${\mathfrak S}^1(H_2,H_1)$ and therefore to 
$\Gamma(X,Y)$.
The duality between $\Gamma(X,Y)$ and 
$B(H_1,H_2)$ is 
given   by $$\langle T,F\rangle=\sum_{n=1}^{\infty}(Tf_n,\bar{g}_n),$$ 
with $T\in B(H_1,H_2)$ and 
$F(x,y)=\sum_{n=1}^{\infty}f_n(x)g_n(y)$.
This will allow us to introduce the
notion of ``operator'' synthesis for some sets of pairs of projections -
bilattices - which  (for separable $H_i$) bijectively correspond to
$\omega$-closed subsets in the product of measure spaces.

Before we proceed with this we give two more definitions which will be used
later.

\begin{definition1}
A synthetic pseudo-closed set is called (operator) 
solvable if each its  pseudo-closed subset is synthetic.
\end{definition1}

Let $X$, $Y$ be standard Borel sets (without measures). We say that 
$K\subseteq X\times Y$ is {\it universally pseudo-closed} if $K$ is 
the complement of a countable union of Borel rectangles. Note that 
if $X$, $Y$ are topological spaces with the natural Borel structure 
then  any closed subset is universally pseudo-closed.

\begin{definition1}
A universally pseudo-closed set $K\subseteq X\times Y$ is said to be 
universally synthetic if
it is $\mu\times\nu$-synthetic for any pair $(\mu,\nu)$ of finite measures.
\end{definition1}

\section{ Bilattices, bimodules and operator synthesis}
First, we introduce the concept of a bilattice and give some notations.   
Let ${\mathcal P}(H)$ denote the lattice of all orthogonal projections
in $B(H)$, the algebra of bounded operators on a Hilbert space $H$. More 
generally, for a von Neumann algebra ${\mathcal R}\subseteq
B(H)$ we denote by ${\mathcal P}_{\mathcal R}$ the lattice of all 
orthogonal projections in ${\mathcal R}$ (thus ${\mathcal P}_{\mathcal R}=
{\mathcal R}\cap{\mathcal P}(H)$, ${\mathcal P}(H)={\mathcal P}_{B(H)}$).

Let $H_1$, $H_2$ be Hilbert spaces. A subset $S\subseteq{\mathcal P}(H_1)\times
{\mathcal P}(H_2)$ is called a {\it bilattice} if
\begin{itemize}
\item $(0,0)$, $(0,1)$, $(1,0)\in S$;
\item $(P_1,Q_1)$, $(P_2,Q_2)\in S\Rightarrow 
(P_1\bigwedge P_2,Q_1\bigvee Q_2)$, $(P_1\bigvee P_2,Q_1\bigwedge Q_2)\in S$.
\end{itemize}
For a bilattice $S$ we denote by $S_l$ and $S_r$ 
the projections of $S$ to 
${\mathcal P}(H_1)$ and ${\mathcal P}(H_2)$ respectively. Clearly,
$S_l$ and $S_r$ are lattices of projections containing $0$ and $1$.
\begin{lemma1}\label{l}
(i) $S_l=\{P\mid (P,0)\in S\}$, $S_r=\{Q\mid (0,Q)\in S\}$.

(ii) If $(P,Q)\in S$, $P_1\leq P$, $Q_1\leq Q$ and $P_1\in S_l$, $Q_1\in S_r$, 
then $(P_1,Q_1)\in S$.
\end{lemma1}  
\begin{proof}
(i) follows from the equality $(P,0)=(P\bigvee 0,Q\bigwedge 0)$.

By (i) $(P_1,0)\in S$, $(0,Q_1)\in S$, whence 
$(P_1,Q)=(P\bigwedge P_1,Q\bigvee 0)\in S$ and $(P_1,Q_1)=(P_1\bigvee 0, 
Q\bigwedge Q_1)\in S$. 
\end{proof}

In what follows we consider only bilattices closed in the strong operator 
topology.  By Lemma~\ref{l},(i), in this case the lattices $S_l$ 
and $S_r$ are
also strongly closed.

To see examples of bilattices note that any subset $U$ of $B(H_1,H_2)$ defines 
a strongly closed bilattice
$$\text{Bil } U=\{(P,Q)\in {\mathcal P}(H_1)\times {\mathcal P}(H_2)\mid 
QTP=0\text{ for any }T\in U\}.$$
Conversely, given a subset ${\mathcal F}\subseteq {\mathcal P}(H_1)\times
{\mathcal P}(H_2)$  we set 
$${\mathfrak M}({\mathcal F})=\{T\in B(H_1,H_2)\mid QTP=0 \text{ for each }
(P,Q)\in{\mathcal F}\}.$$
These maps are in a Galois duality:
$$\text{Bil }({\mathfrak M}(\text{Bil } U))=\text{Bil } U,
\quad {\mathfrak M}(\text{Bil }{\mathfrak M}({\mathcal F}))={\mathfrak M}
({\mathcal F}).$$

It is not difficult to see that  spaces of the form 
${\mathfrak M}({\mathcal F})$ are exactly the reflexive (in sense of
\cite{logsh}) operator spaces; 
 they are characterized by the equality $U={\mathfrak M}(\text{Bil }U)$.
Similarly, the bilattices of the form $\text{Bil } U$ are characterized by
the equality $S=\text{Bil }{\mathfrak M}(S)$ and can be called reflexive.

It is easy to check that ${\mathfrak M}(S)$ is a bimodule over the algebras
${\mathcal A}_l=\text{alg }S_l$, 
${\mathcal A}_r=(\text{alg }S_r)^*$:
$${\mathcal A}_r{\mathfrak M}(S){\mathcal A}_l\subseteq{\mathfrak M}(S).$$

The partial converse of this fact is following:
if $U\subseteq B(H_1,H_2)$ is a bimodule over unital 
subalgebras $W_1\subseteq B(H_1)$, 
$W_2\subseteq B(H_2)$ then any pair $(P,Q)\in \text{Bil }U$ is majorized by 
a pair $(P',Q')\in\text{Bil } U\cap (\text{lat }W_1,\text{lat }W_2^*)$.
Indeed, one sets $P'H=\overline{W_1PH}$, $Q'H=\overline{W_2^*QH}$.

Let ${\mathcal R}_1\subseteq B(H_1)$, ${\mathcal R}_2\subseteq B(H_2)$ be von 
Neumann algebras.  A bilattice $S$ is called 
${\mathcal R}_1\times{\mathcal R}_2$-{\it bilattice}
if $S_l={\mathcal P}_{{\mathcal R}_1}$ and $S_r={\mathcal P}_{{\mathcal R}_2}$.
For example, $\text{Bil }U$ is always a $B(H_1)\times B(H_2)$-bilattice.

The above argument shows that if $S$ is an 
${\mathcal R}_1\times{\mathcal R}_2$-bilattice then ${\mathfrak M}(S)$ is an
${\mathcal R}_1'\times{\mathcal R}_2'$-bimodule. Conversely, for an
${\mathcal R}_1'\times{\mathcal R}_2'$-bimodule $U$ we will consider an 
${\mathcal R}_1\times{\mathcal R}_2$-bilattice
$$\text{Bil}_{{\mathcal R}_1,{\mathcal R}_2}U=(\text{Bil }U)\cap
{\mathcal R}_1\times{\mathcal R}_2.$$
If ${\mathcal R}_1$, ${\mathcal R}_2$ are clear we write $\text{bil }U$ instead
 of $\text{Bil}_{{\mathcal R}_1,{\mathcal R}_2}U$.

We will need a bilattice version of  Arveson's reflexivity theorem for
CSL \cite{arv}. Let us call a bilattice $S$ commutative if 
$S_l$ and $S_r$ are commutative.

\begin{theorem1}\label{tmax}
If $S$ is a commutative bilattice then 
$$(\text{Bil }{\mathfrak M}(S))\cap (S_l\times S_r)=S.$$
\end{theorem1}

\begin{proof} It can be reduced, by a $2\times 2$-matrix 
trick, to Arveson's theorem on reflexivity of commutative subspace lattices,
\cite{arv}
(for a coordinate-free proof see \cite{davidson2} or \cite{Sh}).
Indeed, consider the set,
 ${\mathcal L}$, of all projections $\left(\begin{array}{cc}
P&0\\
0&1-Q
\end{array}\right)\in{\mathcal P}_{B(H_1\oplus H_2)}$, where $(P,Q)\in S$.
Clearly, ${\mathcal L}$ is a commutative strongly closed lattice.
 Therefore, ${\mathcal L}$ is reflexive,
i.e., $\text{lat alg }{\mathcal L}={\mathcal L}$.  
One   easily checks that 
$$\text{alg }{\mathcal L}=\{T=(T_{ij})_{i,j=1}^2\in B(H_1\oplus H_2)\mid 
T_{11}\in\text{alg }S_l, T_{22}\in(\text{alg }S_r)^*,
 T_{21}\in{\mathfrak M}(S),
T_{12}=0\}.$$ 
Therefore, if $Q{\mathfrak M}(S)P=0$ for some $P\in S_l$, 
$Q\in S_r$ then $P\oplus (1-Q)\in 
\text{lat alg }{\mathcal L}={\mathcal L}$,
i.e. $(P,Q)\in S$. This yields $(\text{Bil }{\mathfrak M}(S))
\cap (S_l\times S_r)\subseteq S$.
The reverse inclusion is obvious.
\end{proof}

We have, in particular, $\text{Bil}_{{\mathcal D}_1,{\mathcal D}_2}
{\mathfrak M}(S)=S$ for any ${\mathcal D}_1\times{\mathcal D}_2$-bilattice
$S$, where ${\mathcal D}_i$ are commutative von Neumann algebras. Such 
bilattices will be  the main object of this paper. In what follows we suppose that
${\mathcal D}_1$, ${\mathcal D}_2$ are fixed and $\text{bil }U$ means
$\text{Bil}_{{\mathcal D}_1,{\mathcal D}_2} U$ for $U\subseteq B(H_1,H_2)$.

We see that ${\mathfrak M}(S)$ is the largest among all
${\mathcal D}_1'\times{\mathcal D}_2'$-bimodules $U$ with $\text{bil }U=S$.
Now we are going to present the smallest one.

Given a state $\varphi$ on $B(l_2)$, consider a slice operator
$L_{\varphi}: B(l_2\otimes H_1,l_2\otimes H_2)\rightarrow B(H_1,H_2)\}$ 
defined by $L_{\varphi}(A\otimes B)=\varphi(A)B$. 
Let $\text{conv }S$ denote the convex hull of $S$ (in $B(H_1)\times B(H_2)$),
$\text{Conv }S$ the weak (or uniform, see Lemma~\ref{th3a2}) closure of 
$\text{conv }S$  and let 
${\mathcal R}_1=B(l_2)\bar\otimes {\mathcal D}_1$, 
${\mathcal R}_2=B(l_2)\bar\otimes {\mathcal D}_2$. Set
$$F_S=\{(A,B)\in B(l_2\otimes H_1)\times B(l_2\otimes H_2)\mid
(L_{\varphi}(A), L_{\varphi}(B))\in \text{Conv }S \text{ for any }\varphi\}$$
$$\tilde S=\{(P,Q)\in F_S\mid
P,Q \text{ are projections }\}$$ and
define 
$${\mathfrak M}_0(S)=\{X\in B(H_1, H_2)\mid 1\otimes X\in {\mathfrak M}(\tilde S)\},$$
where $1$ is the identity operator on $l_2$. 
Then $F_S\subseteq{\mathcal R}_1\times{\mathcal R}_2$ by the Fubini property
of tensor product (\cite{takesaki}) and ${\mathfrak M}_0(S)$ is an
ultraweakly closed ${\mathcal D}_1'\times{\mathcal D}_2'$-bimodule like
${\mathfrak M}(S)$.
Here and subsequently $\text{bil } 1\otimes U$  for $U\subseteq B(H_1,H_2)$ 
means
$\text{bil}_{{\mathcal R}_1,{\mathcal R}_2} 1\otimes U$, where
${\mathcal R}_i=B(l_2)\bar\otimes {\mathcal D}_i$.

\begin{lemma1}\label{th3a2}
Let $S$ be a commutative ${\mathcal D}_1\times{\mathcal D}_2$-bilattice.
Then
\begin{gather*}
\overline{\text{\rm conv }S}^u =\overline{\text{\rm conv }S}^w=
\{(A,B)\in {\mathcal D}_1\times
{\mathcal D}_2\mid 0\leq A\leq 1, 0\leq B\leq 1, \\
(E_A([\alpha,1]), 
E_B([\beta,1]))\in S, \alpha+\beta>1\}.
\end{gather*}
where ``u'' and ``w'' indicate the ``uniform'' and the ``weak operator 
topology'' closure of the convex hull, $\text{\rm conv }S$, of $S$ and 
$E_X(\cdot)$ is the spectral projection measure of selfadjoint operator $X$.
\end{lemma1}
\begin{proof}
Let ${\mathfrak R}$ denote the set to the right. To see 
that ${\mathfrak R}\subseteq \overline{\text{conv }S}^u$,
set
$\displaystyle A_n=\sum_{i=1}^n\frac{1}{n}E_A([\frac{i}{n},1])$, 
$\displaystyle B_n=\sum_{i=1}^n\frac{1}{n}E_B([\frac{i}{n},1])$ for 
$(A,B)\in {\mathfrak R}$. Clearly, $A_n\to A$ and
$B_n\to B$ uniformly as $n\to\infty$. 
Then, since
$\displaystyle (E_A([\frac{i}{n},1]), E_B([\frac{n-i+1}{n},1])\in S$ and 
$$\displaystyle(A_n,B_n)=\frac{1}{n}\sum_{i=1}^n((E_A([\frac{i}{n},1]), 
E_B([\frac{n-i+1}{n},1])),$$ we have $(A_n,B_n)\in \text{conv }S$ and therefore
$(A,B)\in  \overline{\text{conv }S}^u$.

Next claim is that ${\mathfrak R}$ is convex. In fact, for $(A_1,B_1)$, 
$(A_2,B_2)\in {\mathfrak R}$, 
we have
\begin{gather}
E_{(A_1+A_2)/2}([\alpha,1])=\bigvee_{n}E_{A_1}([\varepsilon_n,1])E_{A_2}
([2\alpha-\varepsilon_n,1]),\nonumber\\
E_{(B_1+B_2)/2}([\beta,1])=\bigvee_{m}E_{B_1}([\varepsilon_m,1])E_{B_2}
([2\beta-\varepsilon_m,1]),\nonumber
\end{gather}
where $\alpha,\beta\in [0,1)$, $\{\varepsilon_n\}$ is a countable dense 
subset of $[0,1]$.
Fix  $\alpha$, $\beta$ such that $\alpha+\beta>1$. Then for $n$, 
$m\in {\mathbb Z}^+$, we have either $\varepsilon_n+\varepsilon_m>1$ which 
gives $(E_{A_1}([\varepsilon_n,1]),E_{B_1}([\varepsilon_m,1]))\in S$ and 
therefore
$(E_{A_1}([\varepsilon_n,1])E_{A_2}([2\alpha-\varepsilon_n,1]), 
E_{B_1}[\varepsilon_m,1])E_{B_2}([2\beta-\varepsilon_m,1,1]))\in S$, or
$(2\alpha-\varepsilon_n)+(2\beta-\varepsilon_m)>1$ which  implies
$(E_{A_1}([\varepsilon_n,1])E_{A_2}([2\alpha-\varepsilon_n,1])), 
E_{B_1}([\varepsilon_m,1])E_{B_2}([2\beta-\varepsilon_m,1]))\in S$.
Since $S$ is a bilattice, $$(E_{(A_1+A_2)/2}([\alpha,1]),
E_{(B_1+B_2)/2}([\beta,1]))\in S.$$

Next step is to prove that ${\mathfrak R}$ is weakly closed. 
Since it is convex it is enough to prove that it is strongly closed.
Let $\{(A_n, B_n)\}\subset {\mathfrak R}$ be a sequence strongly
converging  to $(A,B)\in{\mathcal D}_1\times{\mathcal D}_2$. Then, for
any 
$\varepsilon>0$ and $\alpha$, $\beta<1$, we have 
$$E_A([\alpha,1])\leq s.\lim_{n\to\infty} E_{A_n}([\alpha+\varepsilon,1]) \text{ and }
E_B([\beta,1))\leq s.\lim_{n\to\infty} E_{B_n}([\beta+\varepsilon,1])$$
(the strong limit). 
Since
$( E_{A_n}([\alpha+\varepsilon,1]), E_{B_n}([\beta+\varepsilon,1])\in S$ if 
$\alpha+\beta>1$ and
$S$ is decreasing and closed in the strong operator topology, we obtain 
$(E_A([\alpha,1]),E_B([\beta,1])\in S$. 
If one of $\alpha$, $\beta$ equals $1$, then that
$(E_A([\alpha,1]),E_B([\beta,1])\in S$ follows from 
$\displaystyle E_A(\{1\})=s.\lim_{\varepsilon\to 0}E_A([1-\varepsilon,1])$, 
$\displaystyle E_B(\{1\})=s.\lim_{\varepsilon\to 0}E_B([1-\varepsilon,1])$.
So we can conclude that $(A,B)\in {\mathfrak R}$.

We have therefore 
$$S\subseteq {\mathfrak R}\subseteq\overline{\text{conv }S}^u \subseteq
\overline{\text{conv }S}^w,$$
and, since ${\mathfrak R}$ is convex and weakly closed,  ${\mathfrak R}
=\overline{\text{conv }S}^u=\overline{\text{conv }S}^w$.
\end{proof}

\begin{definition1}
We say that a ${\mathcal D}_1\times{\mathcal D}_2$-bilattice, $S$, is
{\it synthetic} if there exists only one ultraweakly closed 
${\mathcal D}_1'\times{\mathcal D}_2'$-bimodule $\mathfrak M$ such that
$\text{bil }{\mathfrak M}=S$. 
\end{definition1}

\begin{theorem1}\label{th3b} 
Let $S$ be a ${\mathcal D}_1\times{\mathcal D}_2$-bilattice and  
let ${\mathfrak M}$ be an ultraweakly closed 
${\mathcal D}_1'\times{\mathcal D}_2'$-bimodule
such that $\text{\rm bil }{\mathfrak M}\subseteq S$. 
Then $\text{\rm bil }1\otimes{\mathfrak M}\subseteq\tilde S$. 
\end{theorem1}
\begin{proof}
Let $(P,Q)\in \text{bil }1\otimes{\mathfrak M}$. Fix $\xi\in l_2$, $||\xi||=1$.
Consider the corresponding state $\varphi_{\xi}(A)=(A\xi,\xi)$ and denote
the corresponding operator $L_{\varphi_{\xi}}$ simply by $L_{\xi}$.
It suffices to show that $(L_{\xi}(P), L_{\xi}(Q))\in \text{Conv }S$.
By definition of  $L_{\xi}$, we have 
$(L_{\xi}(K)x,x)=(K(\xi\otimes x),\xi\otimes x)$ for any operator $K$ on 
$l_2\otimes H$ and,
in particular, 
if $K=P$ (a selfadjoint projection) then
$(L_{\xi}(P)x,x)=||P(\xi\otimes x)||^2$. Therefore, for $A\in{\mathfrak M}$ the following holds
\begin{eqnarray*}
&(AL_{\xi}(P)A^*x,x)=(L_{\xi}(P)A^*x,A^*x)=||P(\xi\otimes A^*x)||^2= \\
&||P(1\otimes A^*)Q^{\perp}(\xi\otimes x)||^2\leq ||A||^2||Q^{\perp}(\xi\otimes x)||^2=
||A||^2(L_{\xi}(Q^{\perp})x,x). 
\end{eqnarray*}
We obtain now the inequality
$AL_{\xi}(P)A^*\leq ||A||^2L_{\xi}(Q^{\perp})$. Let $L_{\xi}(P)=K^2$, $L_{\xi}(Q^{\perp})=L^2$, where
$K$, $L\geq 0$. Then $||KA^*x||\leq||A||||Lx||$ for any $A\in{\mathfrak M}$ and $x\in H$. If $L$ is invertible this is equivalent to 
$||KA^*L^{-1}||\leq ||A||$.
Since ${\mathfrak M}$ is a bimodule, $KA^*L^{-1}\in{\mathfrak M}^*$. Writing
now
$KA^*L^{-1}$ instead of $A^*$ we get 
$||K^2A^*L^{-2}||\leq||KA^*L^{-1}||\leq ||A||$.
Proceeding in this fashion we obtain
$||K^nA^*L^{-n}||\leq||A||$   and  hence 
\begin{equation}\label{eqq}
||K^nA^*x||\leq||A|| ||L^nx||,\quad x\in H.
\end{equation}
If $L$ is not invertible, then replacing $L$ by $L+\e 1$ in the above argument
we obtain (\ref{eqq}) for all $L+\e 1$ with $\e>0$. Letting $\e\to 0$,
we get (\ref{eqq}) for $L$.

Fix $x\in E_{L}([0,\varepsilon])$, where $E_L(\cdot)$ is the spectral 
projection measure of $L$.
Then $||L^nx||\leq C\varepsilon^n$ and, by (\ref{eqq}), we obtain
 $A^*x\in E_K([0,\varepsilon])$. 
Thus
${\mathfrak M}^* E_{L}([0,\varepsilon])\subseteq  E_K([0,\varepsilon])$
or, equivalently,
$E_K([\varepsilon', 1]){\mathfrak M}^* E_{L}([0,\varepsilon])=0$  if  
$\varepsilon'>\varepsilon$.
This implies $E_{K^2}([\varepsilon',1]){\mathfrak M}^*
E_{1-L^2}([1-\varepsilon,1])=0$, as
 $\varepsilon'>\varepsilon$, i.e. 
$$E_{L_{\xi}(Q)}([\alpha,1]){\mathfrak M}E_{L_{\xi}(P)}([\beta,1])=0, \quad 
\alpha+\beta>1.$$
Since $(E_{L_{\xi}(Q)}([\alpha,1]),E_{L_{\xi}(P)}([\beta,1]))\in 
\text{bil }{\mathfrak M}\subseteq S$ as
$\alpha+\beta>1$, by Lemma~\ref{th3a2}, we obtain
 $(L_{\xi}(Q), L_{\xi}(P))\in \text{Conv }S$ for any 
$\xi\in H$. 
\end{proof}
The idea of the proof goes back to Arveson \cite{arv}.
\begin{cor1}\label{cor0b}
 Let $S$ be a ${\mathcal D}_1\times{\mathcal D}_2$-bilattice and  
let ${\mathfrak M}$ be an ultraweakly closed 
${\mathcal D}_1'\times{\mathcal D}_2'$-bimodule
such that $\text{\rm bil }{\mathfrak M}\subseteq S$. Then 
${\mathfrak M}_0(S)\subseteq{\mathfrak M}$.
\end{cor1}
\begin{proof}
Let $T\in{\mathfrak M}_0(S)$. To see that $T\in{\mathfrak M}$ we choose
an ultraweakly continuous linear functional $\varphi$ such that
$\varphi({\mathfrak M})=0$. Then there exist $F\in l_2\otimes H_1$, $G\in l_2
\otimes H_2$ such that $\varphi(A)=((1\otimes A)F,G)$, 
$A\in B(H_1, H_2)$, moreover, $(1\otimes{\mathfrak M})F\perp G$. Denoting
by $P_F$ and $P_G$ the projections on $\overline{[(1\otimes{\mathcal D}_1')F]}$
 and $\overline{[(1\otimes{\mathcal D}_2')G]}$
we have $P_G(1\otimes{\mathfrak M})P_F=0$, i.e. \ $(P_F,P_G)\in 
\text{bil }1\otimes{\mathfrak M}$. It follows now from the definition of 
${\mathfrak M}_0(S)$ and Theorem~\ref{th3b} that 
$(P_F,P_G)\in \tilde S\subseteq 
\text{bil }1\otimes T$ and therefore $P_G(1\otimes T)P_F=0$, i.e. $\varphi(T)=0$.
 From the arbitrariness
 of $\varphi$ we obtain $T\in{\mathfrak M}$. 
\end{proof}
Summarising we have the  following statement. 
\begin{theorem1}\label{thad}
Let $S$ be a ${\mathcal D}_1\times{\mathcal D}_2$-bilattice. If 
${\mathfrak M}$ is an ultraweakly closed 
${\mathcal D}_1'\times{\mathcal D}_2'$-bimodule
such that $\text{\rm bil }{\mathfrak M}=S$ then 
${\mathfrak M}_0(S)\subseteq{\mathfrak M}\subseteq{\mathfrak M}(S)$.
\end{theorem1}

\begin{theorem1}\label{smallest} 
Given  a ${\mathcal D}_1\times{\mathcal D}_2$-bilattice 
$S$, 
$\text{bil }{\mathfrak M}_0(S)=S$.
\end{theorem1}

Theorem~\ref{thad} and \ref{smallest} state that
${\mathfrak M}_0(S)$ is the smallest ultraweakly closed 
${\mathcal D}_1'\times{\mathcal D}_2'$-bimodule whose bilattice is $S$ 
 and that a commutative 
bilattice $S$ is
synthetic if and only if ${\mathfrak M}(S)={\mathfrak M}_0(S)$.

We shall prove Theorem~\ref{smallest} in Section~\ref{sssec122} after treating
the case of  bilattices on separable Hilbert spaces.
Here we only give one of its consequences.
\begin{cor1}
If ${\mathcal L}$ is a CSL then there is a smallest element in the class of all
ultra-weakly closed algebras ${\mathcal A}$ such that 
$\text{lat }{\mathcal A}={\mathcal L}$ and ${\mathcal L}'\subseteq{\mathcal A}$.
\end{cor1}
\begin{proof}
Set ${\mathcal D}={\mathcal L}''$ and 
$$S=\{(P,Q)\in{\mathcal P}_{\mathcal D}\times{\mathcal P}_{\mathcal D}\mid
\exists R\in{\mathcal L}\text{ with } P\leq R\leq 1-Q\}.$$
Then $S$ is a ${\mathcal D}\times{\mathcal D}$-bilattice. We denote by
${\mathcal A}_0({\mathcal L})$ the ultra-weakly closed algebra generated 
by ${\mathfrak M}_0(S)$. 

Note that $1\in{\mathfrak M}_0(S)$. Indeed, since
$P+Q\leq 1$ for any $(P,Q)\in S$, $$\text{Conv }S\subseteq\{(A,B)\in
{\mathcal D}\times{\mathcal D}\mid A+B\leq 1\}.$$ Hence if $(P_1,Q_1)\in\tilde S$
 then $L_{\varphi}(P_1+Q_1)\leq 1$ for any state $\varphi$ on $B(l_2)$.
Since $P_1+Q_1\in B(l_2)\bar\otimes{\mathcal D}$,  
we can conclude, using \cite{EKS}[Lemma~7.5, (ii)], that
 $P_1+Q_1\leq 1$. Hence $P_1Q_1=0$ and $Q_1(1\otimes 1)P_1=0$,
$1\in {\mathfrak M}_0(S)$.

Since ${\mathfrak M}_0(S)$ is an ${\mathcal L}'$-bimodule, 
we have ${\mathcal L}'\subseteq{\mathfrak M}_0(S)$ and 
${\mathcal L}'\subseteq{\mathcal A}_0({\mathcal L})$.
Let us show that $\text{lat }{\mathcal A}_0({\mathcal L})={\mathcal L}$.
Indeed, $$\text{lat }{\mathcal A}_0({\mathcal L})\subseteq\text{lat }
({\mathcal L}')\subseteq{\mathcal D}$$
 and $P\in{\mathcal P}_{\mathcal D}$ belongs to $\text{lat }{\mathcal A}_0
({\mathcal L})$ iff $P\in\text{lat }{\mathfrak M}_0(S)$ iff
$(1-P){\mathfrak M}_0(S)P=0$ iff $(P,1-P)\in\text{bil }{\mathfrak M}_0(S)=S$
iff $P\leq R\leq P$ for some $R\in{\mathcal L}$ iff $P\in{\mathcal L}$.

Let ${\mathcal A}$ be an ultra-weakly closed algebra containing ${\mathcal D}$
and $\text{lat }{\mathcal A}={\mathcal L}$. Then $\text{bil }{\mathcal A}=S$.
Indeed, if $(P,Q)\in S$ then there is $R\in{\mathcal L}$ such that
$P\leq R$, $Q\leq 1-R$, whence
$$Q{\mathcal A}P=Q(1-R){\mathcal A}RP=0,$$
$(P,Q)\in \text{bil }{\mathcal A}$. Conversely, if
$(P,Q)\in\text{bil }{\mathcal A}$ then
setting $RH=\overline{{\mathcal A}PH}$ we have $R\in{\mathcal L}$,
$Q{\mathcal A}R=0$ whence $QR=0$, $Q\leq 1-R$ and 
$(P,Q)\in S$, because $P\leq R$.

By Theorem~\ref{thad}, ${\mathfrak M}_0(S)\subseteq{\mathcal A}$ whence
${\mathcal A}_0({\mathcal L})\subseteq{\mathcal A}$. 
\end{proof}
\begin{remark}\rm
Arveson \cite{arv} calls an ultra-weakly closed algebra ${\mathcal A}$
with $\text{lat }{\mathcal A}={\mathcal L}$ {\it pre-reflexive} if
${\mathcal L}'\subseteq{\mathcal A}$. In this terms corollary can
be considered as an extension to non-separable spaces of the result
by Arveson \cite{arv}[Theorem~2.1.8, (ii)] on the existence of the
smallest pre-reflexive algebra with a given commutative lattice.
\end{remark}

\section{ Separably acting  bilattices}\label{sec4}
If Hilbert spaces $H_1$ and $H_2$ are separable then there exist finite
separable measure spaces $(X,\mu)$ and $(Y,\nu)$ with  standard measures $\mu$, $\nu$,
such that $H_1=L_2(X,\mu)$, $H_2=L_2(Y,\nu)$
and the multiplication algebras ${\mathcal D}_1$ and ${\mathcal D}_2$ 
are $L_{\infty}(X,\mu)$ and $L_{\infty}(Y,\nu)$ 
respectively. Denote by
$P_U$ and $Q_V$   the multiplication operators by 
the characteristic functions of $U\subseteq X$ and $V\subseteq Y$.
Given $E\subseteq X\times Y$, we define  $S_E$ to be the set of all
pairs of projections $(P_U,Q_V)$, where $U\subseteq X$, $V\subseteq Y$ and
$(U\times V)\cap E\cong\emptyset$. 

\begin{theorem1}
 $S_E$ is a ${\mathcal D}_1\times{\mathcal D}_2$-bilattice. 
\end{theorem1}
\begin{proof}
We shall prove only the closedness of $S_E$, the other conditions 
trivially hold. Let $(P_n,Q_n)\in S_E$, $P_n\to P$, $Q_n\to Q$ in the strong 
operator topology. Then there exist $A\subseteq X$, $B\subseteq Y$ such that
$P=P_A$, $Q=Q_B$. Changing, if necessarily, $P_n$ to $P_nP$, $Q_n$ to
$Q_nQ$, we may assume that $P_n\leq P$ and $Q_n\leq Q$. We have therefore
$P_n=P_{A_n}$, $Q_n=Q_{B_n}$, for some $A_n\subseteq X$, $B_n\subseteq Y$
such that $(A_n\times B_n)\cap E\cong\emptyset$ and 
$\mu(A\setminus A_n)\to 0$,  $\nu(B\setminus B_n)\to 0$.
Given $\e>0$, $k\in{\mathbb N}$, choose $n_k$ such that
$\displaystyle\mu(A\setminus A_{n_k})<\frac{\e}{2^k}$ and 
$\displaystyle\nu(B\setminus B_{n_k})<\frac{\e}{2^k}$. Set 
$$\displaystyle A_{\e}=\cap_{k=1}^{\infty}A_{n_k}, \quad 
B_{\e}=\cup_{k=1}^{\infty}B_{n_k}.$$
Then $\mu(A\setminus A_{\e})\leq \e$, $\nu(B\setminus B_{\e})=0$ and 
$(A_{\e}\times B_{\e})\cap E\cong\emptyset$. Taking now
$A_0=\cup_{n=1}^{\infty}A_{1/n}$ and $B_0=\cap_{n=1}^{\infty}B_{1/n}$, we
obtain
$\mu(A\setminus A_0)=0$, $\nu(B\setminus B_0)=0$, $(A_0\times B_0)\cap E\cong 
\emptyset$ so that $(P,Q)=(P_{A_0}, Q_{B_0})\in S_E$.

\end{proof}

\begin{theorem1}
Let  $S$ be a ${\mathcal D}_1\times{\mathcal D}_2$-bilattice. Then
 there exists a unique, up to a marginally null set,
 pseudo-closed set $E\subseteq X\times Y$  such that $S=S_E$.
\end{theorem1}
\begin{proof}
Let $\{(P_n,Q_n)\}$ be a strongly dense sequence in the bilattice $S$, and let
$A_n\subseteq X$, $B_n\subseteq Y$ be such that $P_n=P_{A_n}$ and 
$Q_n=Q_{B_n}$. The set $\displaystyle E=(X\times Y)\setminus 
(\cup_{n=1}^{\infty}A_n\times B_n)$ is clearly pseudo-closed. We will
show that $S=S_{E}$. 

Since $S_E$ is closed in the strong operator topology,
we have the inclusion $S\subseteq S_E$. For the reverse inclusion, we first
show  that if a  rectangle, $A\times B$, lies in the union of a finite
number of rectangles, say $C_k\times D_k$ ($1\leq k\leq n$), such that 
$(P_{C_k},Q_{D_k})\in S$, then $(P_A,Q_B)\in S$.
We use the induction by $n$. The case $n=1$ is obvious from the 
decreasing condition on $S$. If 
$A\times B\subseteq\cup_{k=1}^n C_k\times D_k$,
then $(A\setminus C_1)\times B\subseteq\cup_{k=2}^n (C_k\times D_k)$ 
and so, by the induction hypothesis, we have that
$(P_{A\setminus C_1},Q_B)\in S$. Similarly, $(P_A,Q_{B\setminus D_1})\in S$.
Therefore, $(P_{A\cap C_1},Q_{B\setminus D_1})\in S$. Since
$S$ is closed under the operation $(\bigvee,\bigwedge)$, this   together
with $(P_{C_1},P_{D_1})\in S$ gives us $(P_{A\cap C_1},P_B)\in S$.
Using again  closeness under $(\bigwedge,\bigvee)$, we
obtain $(P_A,Q_B)\in S$.

Let now $(P,Q)=(P_A,Q_B)\in S_E$.  Deleting null sets from $A$, $B$ 
we may assume that $A\times B\subseteq\cup_{n=1}^{\infty}
A_n\times B_n$. Then, by \cite{EKS}[Lemma~3.4,d], given $\e>0$, there
exist $A_{\e}\subseteq A$, $B_{\e}\subseteq B$ with
$\mu(A\setminus A_{\e})<\e$, $\nu(B\setminus B_{\e})<\e$ such that
$A_{\e}\times B_{\e}$ is contained in the union of a finite number
of sets $\{A_n\times B_n\}$. By the statement we have just proved, 
$(P_{A_{\e}},Q_{B_{\e}})\in S$, and, since $P_{A_{\e}}\to P$, 
$Q_{B_{\e}}\to Q$ strongly, as $\e\to 0$, we have $(P,Q)\in S$.  
This proves $S=S_E$. 

To see the uniqueness, let $E_1$ be a pseudo-closed set such that 
$S_{E_1}=S_E$. Then $(P_A,Q_B)\in S_E$ for any $A\times B\in E_1^c$ and 
therefore $A\times B\subseteq E^c$ up to a marginally null set.
As $E_1^c$ is pseudo-open, we have $E_1^c\subseteq E^c$ up to a marginally null
set. Similarly, we have the reverse inclusion and
therefore $E_1^c\cong E^c$ and  $E_1\cong E$.
\end{proof}

We say that $T\in B(H_1,H_2)$ is {\it supported} in $E\subseteq 
X\times Y$ if
$\text{bil }T\supseteq S_E$, i.e., if $Q_VTP_U=0$ for each
sets $U\subseteq X$,  $V\subseteq Y$ such that $(U\times V)\cap E\cong
\emptyset$. Clearly, 
$${\mathfrak M}(S_E)=\{T\in B(H_1,H_2)\mid\text{$T$ is supported in } E\}.$$
For any  subset ${\mathbb U}\subseteq B(H_1,H_2)$ there exists the smallest
(up to a marginally null set ) pseudo-closed set, $\supp {\mathbb U}$, 
which supports any operator $T\in {\mathbb U}$, namely, $\supp {\mathbb U}$ 
is the pseudo-closed set $E$ such that $\text{bil }{\mathbb U}=S_E$.
The support of an operator $T\in B(H_1,H_2)$ will be  denoted by $\supp T$. 
  We will also use  the notations
${\mathfrak M}_{max}(E)$ and ${\mathfrak M}_{min}(E)$ for the 
bimodules ${\mathfrak M}(S_E)$ and  ${\mathfrak M}_0(S_E)$.
Theorem~\ref{thad} says now that
$${\mathfrak M}_{min}(E)\subseteq{\mathfrak M}\subseteq{\mathfrak M}_{max}(E)$$
if $\supp {\mathfrak M}=E$. Clearly,
$\supp {\mathfrak M}_{max}(E)=E$ and therefore ${\mathfrak M}_{max}(E)$ is 
the largest  ultraweakly closed bimodules whose support is $E$. By proving 
now that $\supp {\mathfrak M}_{min}(E)=E$ we would also have that 
${\mathfrak M}_{min}(E)$ is the smallest ultraweakly closed bimodules whose 
support is $E$, justifying the notations.

Let $\Psi$ be a subspace of $\Gamma(X,Y)$. Using the duality of
$B(H_1,H_2)$ and $\Gamma(X,Y)$ we denote by
 $\Psi^{\perp}$ the subspace of
all operators $T\in B(H_1,H_2)$ such that $\langle T, F\rangle=0$ for any
$F\in\Psi$. Clearly, if $\Psi$ is invariant then $\Psi^{\perp}$ is 
a $({\mathcal D}_1,{\mathcal D}_2)$-bimodule.  

\begin{theorem1}\label{th5}
 Let $E\subseteq X\times Y$ be a pseudo-closed set. Then
$$\Phi_0(E)^{\perp}={\mathfrak M}_{max}(E).$$
\end{theorem1}
\begin{proof}
We begin by showing the inclusion 
${\mathfrak M}_{max}(E)\subseteq\Phi_0(E)^{\perp}$. Let 
$A\in{\mathfrak M}_{max}(E)$, $F\in\Phi_0(E)$. By \cite{EKS}[Lemma~3.4], 
$E$ is $\e$-compact, so that, for any $\e>0$, there exist $X_{\e}\subseteq X$, 
$Y_{\e}\subseteq Y$ with  $\mu(X_{\e})<\e$, $\nu(Y_{\e})<\e$ such that 
$$F_{\e}(x,y)=F(x,y)\chi_{X_{\e}^c}(x)\chi_{Y_{\e}^c}(y)$$ 
vanishes on an open-closed neighbourhood of $E$ ($\cong$ the union of a 
finite number of rectangles). Clearly,
$F_{\e}\to F$ as $\e\to 0$. It remains to show that 
$\langle A, F_{\e}\rangle=0$. Choose measurable sets $\{X_j\}_{j=1}^N$, 
$\{Y_i\}_{i=1}^M$
in a way that $$X=\cup_{j=1}^NX_j,\quad Y=\cup_{i=1}^M Y_i \quad\text{and}
\quad\text{null }F_{\e}\supseteq \cup_{(i,j)\in J}X_j\times Y_i\supseteq E$$ 
for some index set $J$. 
If $(i,j)\in J$ then $\langle Q_{Y_i}AP_{X_j}, F_{\e}\rangle=
\langle A, F_{\e}\chi_{X_j}\chi_{Y_i}\rangle=0$.
If $(i,j)\notin J$ then  $Q_{Y_i}AP_{X_j}=0$ since $\supp A\subseteq E$.
Therefore, $\langle Q_{Y_i}AP_{X_j}, F_{\e}\rangle=0$ for any pair $(i,j)$ 
and hence $\langle A, F_{\e}\rangle=0$.

Let $A$ be an operator in $B(H_1,H_2)$  such that $\langle A,F\rangle=0$
for any $F\in\Phi_0(E)$. Consider $U\subseteq X$, $V\subseteq Y$ such that 
$(U\times V)\cap E=\emptyset$ (up to a marginally null set). Then 
$F(x,y)\chi_U(x)\chi_V(y)\in\Phi_0(E)$
for any $F\in \Gamma(X,Y)$ and
$\langle Q_VAP_U,F\rangle=\langle A, F\cdot\chi_V\chi_U\rangle=0$, 
which implies $Q_VAP_U=0$. 
\end{proof}

Let ${\mathcal D}_i^+$ denote the set of positive functions in 
${\mathcal D}_i$. 
Operators $A\in B(l_2)\bar\otimes{\mathcal D}_1$ and 
$B\in B(l_2)\bar\otimes{\mathcal D}_2$ can be
identified with operator-valued functions $A(x):X\to B(l_2)$ and 
$B(y):Y\to B(l_2)$. 
If $A$, $B$ are projections then $A(x)$, $B(y)$ are projection-valued 
functions. We say that  a pair of projections 
$(P, Q)\in (B(l_2)\bar\otimes{\mathcal D}_1)\times  
(B(l_2)\bar\otimes{\mathcal D}_2)$ is an {\it $E$-pair} if $P(x)Q(y)$ 
vanishes on $E$. If, additionally, $P$ and $Q$ take only finitely many 
values then the pair $(P,Q)$ is said to be a {\it simple $E$-pair}. 

\begin{lemma1} \label{cor0a}
Let $E$ be a pseudo-closed subset of $X\times Y$. Then
$$\text{\rm Conv }S_E=\{(a(x),b(y))\in 
{\mathcal D}_1^+\times{\mathcal D}_2^+\mid
a(x)+b(y)\leq 1, \text{ m.a.e on } E\},$$
$$F_{S_E}=\{(A,B)\in (B(l_2)\bar\otimes{\mathcal D}_1)^+\times  
(B(l_2)\bar\otimes{\mathcal D}_2)^+\mid A(x)+B(y)\leq 1, 
\text{ m.a.e on } E\},$$  
and
$$\tilde S_E=\{(P,Q)\mid (P,Q)\text{ is an $E$-pair}\}.$$ 
\end{lemma1}
\begin{proof}
The first statement follows easily from Lemma~\ref{th3a2}.
To see the second equality take $\xi\in l_2$ and $(A,B)\in F_{S_E}$, 
identifying the operators with the corresponding operator-valued 
functions. 
Set now $a(x)=(A(x)\xi,\xi)$ and
$b(y)=(B(y)\xi,\xi)$. It is easy to see that $(L_{\xi}(A)f)(x)=a(x)f(x)$
and $(L_{\xi}(B)g)(y)=b(y)g(y)$. By the definition of $F_{S_E}$ and 
the first statement, we have 
$(A(x)+B(y)\xi,\xi)=(A(x)\xi,\xi)+(B(y)\xi,\xi)=
a(x)+b(y)\leq 1$ (m.a.e.) on $E$ and therefore $A(x)+B(y)\leq 1$
(m.a.e.)\ on $E$. If, additionally,  $A$ and $B$ are projections, 
the inequality gives $A(x)B(y)=0$ (m.a.e.)\ on $E$, completing
the proof. 
\end{proof}

\begin{theorem1} \label{th6}
Let $E\subseteq X\times Y$ be a pseudo-closed set. Then
$$\Phi(E)^{\perp}={\mathfrak M}_{min}(E).$$
\end{theorem1}
\begin{proof}
Let $(P,Q)\in \tilde S_E$ and let $\vec{x}(x)=P(x)\xi$ and 
$\vec{y}(y)=Q(y)\eta$ for some $\xi$, $\eta\in l_2$. By Lemma~\ref{cor0a},
$(P(x),Q(y))$  is an $E$-pair which implies   
$(\vec{x}(x),\vec{y}(y))=0$ m.a.e.\ on $E$. Clearly, the function 
$F:(x,y)\mapsto (\vec{x}(x),\vec{y}(y))$ belongs to $\Gamma(X,Y)$ and therefore
$F\in\Phi(E)$. For any $T\in B(H_1,H_2)$ we have $\langle T,F\rangle=
((1\otimes T)\vec{x},\vec{y})$ and if $T\in\Phi(E)^{\perp}$ we obtain
$((1\otimes T)\vec{x},\vec{y})=0$ and $Q(1\otimes T)P=0$, i.e. 
$T\in{\mathfrak M}_{min}(E)$.

To see the converse we observe that any function $F\in \Phi(E)$ can be written
as $(\vec{x}(x),\vec{y}(y))$, where $\vec{x}(x)$, $\vec{y}(y)\in l_2$ and
$\vec{x}(x)\perp\vec{y}(y)$ if $(x,y)\in E$ m.a.e. Denoting
by $P(x)$ and $Q(y)$ the projections onto the one-dimensional spaces generated
by $\vec{x}(x)$ and $\vec{y}(y)$ yields $P(x)Q(y)=0$ m.a.e.\ on $E$ and
$(P,Q)\in \tilde S_E$. For any $T\in {\mathfrak M}_{min}(E)$ we have
$$\langle T,F\rangle=((1\otimes T)\vec{x}(x),\vec{y}(y))=
(Q(1\otimes T)P\vec{x}(x),\vec{y}(y))=0.$$ 
This  implies $T\in\Phi(E)$.
\end{proof}
\begin{cor1}\label{cor1b}
$$\text{\rm bil }{\mathfrak M}_{min}(E)=S_E.$$
\end{cor1}
\begin{proof} It suffices to show that
 $Q_V{\mathfrak M}_{min}(E)P_U=0$ with measurable $U\subseteq X$, 
$V\subseteq Y$ implies that $(U\times V)\cap E$ is marginally null. 
In fact, this would imply $S_E\supseteq \text{bil }{\mathfrak M}_{min}(E)$ 
which together 
with $S_E=\text{bil }{\mathfrak M}_{max}(E)\subseteq 
\text{bil }{\mathfrak M}_{min}(E)$ gives
us the statement. The last inclusion holds since 
${\mathfrak M}_{min}(E)\subseteq {\mathfrak M}_{max}(E)$.  

Assume that $E_0=(U\times V)\cap E$ is not marginally null. Then
$\Phi(E_0)$ does not contain $\chi_{U\times V}$ and
therefore is not equal to $\Gamma(U,V)$. Since $\Phi(E_0)$ is closed
in $\Gamma(U,V)$, there exists an operator
$A_0\in B(P_UH_1,Q_VH_2)$ such that $0\ne A_0\perp\Phi(E_0)$. 
Extend $A_0$ to an operator $A\in B(H_1,H_2)$ so that $Q_VAP_U|_{L_2(U)}=A_0$
and $A=Q_VAP_U$. Then $A\perp\Phi(E)$ and, by Theorem~\ref{th6},
$A\in{\mathfrak M}_{min}(E)$.
Since $Q_VAP_U\ne 0$, we obtain a contradiction.
\end{proof}
\begin{cor1}\label{cc1} 
Let ${\mathfrak M}\subseteq B(H_1,H_2)$ be an ultraweakly closed 
bimodule,  $E$ be a pseudo-closed set. 
Then $\supp {\mathfrak M}=E$ iff
$${\mathfrak M}_{min}(E)\subseteq{\mathfrak M}\subseteq
{\mathfrak M}_{max}(E).$$
\end{cor1}
\begin{proof}
It follows from Theorem~\ref{tmax}, Corollary~\ref{cor0b},\ref{cor1b} and the 
fact that
$\text{bil }{\mathfrak M}=S_E$ if and only if $\supp {\mathfrak M}=E$.
\end{proof}
{\bf Proof of Theorem~\ref{th1}.}
 Let $E=\supp A^{\perp}$. By Corollary~\ref{cc1},
$${\mathfrak M}_{min}(E)\subseteq A^{\perp}\subseteq{\mathfrak M}_{max}(E)$$
and therefore, by Theorem~\ref{th5}, \ref{th6},
$$\Phi_0(E)\subseteq A\subseteq\Phi(E)$$
which also implies $\text{null }A=E$. \qed

The next corollary is an analogue of Wiener's Tauberian Theorem.
\begin{cor1}\label{dense}
If $\Psi\subseteq\Gamma(X,Y)$ and $\text{\rm null }\Psi\cong\emptyset$  then 
$\Psi$ is dense in $\Gamma(X,Y)$.
\end{cor1}
\begin{proof}
Follows from Theorem~\ref{th1}, since $\Phi_0(\emptyset)=\Gamma(X,Y)$.
\end{proof}
\begin{cor1}\label{cor1}
$$\text{\rm bil }1\otimes{\mathfrak M}_{min}(E)=\tilde S_E=
\{(P,Q):(P,Q)\text{ is an 
$E$-pair }\},$$
\end{cor1}
\begin{proof} 
By Corollary~\ref{cor1b}, $\text{bil }{\mathfrak M}_{min}(E)=S_E$ which 
together with Theorem~\ref{th3b} implies
$\text{bil }1\otimes{\mathfrak M}_{min}(E)\subseteq\tilde S_E$. 
On the other hand,
$\text{bil }1\otimes{\mathfrak M}_{min}(E)\supseteq\tilde S_E$ by the 
definition of ${\mathfrak M}_{min}(E)$. The second equality is proved in 
Lemma~\ref{cor0a}.
\end{proof}
\begin{remark}\rm For sets that are graphs of preoders (that is for lattices)
the result was, in fact, proved in
\cite{arv}[Cor.1 of Theorem~2.1.5].
\end{remark}

\begin{theorem1}\label{th7}
 Let $E$ be a pseudo-closed set. Then 
$$\text{\rm bil }1\otimes{\mathfrak M}_{max}(E)=
\overline{\{(P,Q):(P,Q)\text{ is a simple $E$ 
pair }\}}^s,$$
where ``s'' indicates the strong operator topology closure.
\end{theorem1}
\begin{proof}
Consider the commutative lattice, ${\mathcal L}$, of all projections
$\left(\begin{array}{cc}
p&0\\
0&1-q
\end{array}\right)\in {\mathcal P}_{B(H_1\oplus H_2)}$, where
$(p,q)\in S_E$. By \cite{Sh}, 
\begin{equation}\label{lat}
{\mathcal P}_{B(l_2)}\otimes{\mathcal L}=\text{lat }(1\otimes
\text{alg }{\mathcal L}),
\end{equation}
where  the tensor product on the left hand side denotes the smallest 
(strongly closed) lattice 
containing the elementary tensors $A\otimes B$, $A\in {\mathcal P}_{B(l_2)}$,
$B\in {\mathcal L}$. 
Moreover, it is shown in \cite{Sh} that
$$\text{lat }(1\otimes
\text{alg }{\mathcal L})=\lim_n {\mathcal P}_{B(l_2)}\otimes{\mathcal L}_n,$$
where ${\{\mathcal L}_n\}$ is a sequence of finite sublattices of 
${\mathcal L}$. It is easy to check that for a finite sublattice 
${\mathcal L}_n\subseteq{\mathcal L}$, 
${\mathcal P}_{B(l_2)}\otimes{\mathcal L}_n\subseteq\{P\oplus (1-Q): (P,Q) \text{ is a 
simple $E$-pair}\}$, whence

$${\mathcal P}_{B(l_2)}\otimes{\mathcal L}=
\overline{\{P\oplus (1-Q): (P,Q) \text{ is a 
simple $E$-pair}\}}^s.$$
 
Since $$\text{alg }{\mathcal L}=\{T=(T_{ij})_{i,j=1}^2\in B(H_1\oplus H_2)\mid
T_{11}\in{\mathcal D}_1,T_{22}\in{\mathcal D}_2,T_{21}\in{\mathfrak M}_{max}(E),
T_{12}=0\}$$ 
(see the proof of Theorem~\ref{tmax}), one can easily check
that $\left(\begin{array}{cc}
P&0\\
0&1-Q
\end{array}\right)\in {\mathcal P}_{B(l_2\otimes H_1\oplus l_2\otimes H_2)}$,
where $(P,Q)\in \text{bil }(1\otimes {\mathfrak M}_{max}(E))$,
belongs to $\text{lat }(1\otimes
\text{alg }{\mathcal L})$.  By (\ref{lat}) we have
$\text {\rm bil }1\otimes{\mathfrak M}_{max}(E)\subseteq
\overline{\{(P,Q):(P,Q)\text{ is a simple $E$ 
pair }\}}^s$. The reverse inclusion is obvious.
\end{proof}
In the following theorem we list several possible definitions of a set 
of operator synthesis.

\begin{theorem1}\label{main}
 Let $E\subseteq X\times Y$ be a pseudo-closed set.
Then the following are equivalent:

$(i)$ $E$ is a set of synthesis;

$(ii)$ ${\mathfrak M}_{min}(E)={\mathfrak M}_{max}(E)$;

$(iii)$ $\langle T,F\rangle=0$ for any $T\in B(H_1,H_2)$ and 
$F\in \Gamma(X,Y)$,
 $\text{\rm \supp} T\subseteq E\subseteq \text{\rm null } F$;

$(iv)$ any $E$-pair can be approximated in the strong operator topology
of $B(l_2\otimes H_1)\times B(l_2\otimes H_2)$ by simple $E$-pairs;

$(v)$ any $E$-pair can be approximated by simple $E$-pairs
almost everywhere in the strong operator topology of $B(l_2)$. 
\end{theorem1}
\begin{proof}
$(i)\Leftrightarrow (ii)$: obviously follows from the definition and
Theorems~\ref{th5},\ref{th6}.

$(ii)\Rightarrow (iii)$: if $T\in{\mathfrak M}_{min}(E)$ then, by 
Theorem~\ref{th6}, $\langle T,F\rangle=0$ for any $F\in\Gamma(X,Y)$,
such that $E\subseteq \text{null }F$, which shows the implication.

$(iii)\Rightarrow (ii)$: Let $T\in {\mathfrak M}_{max}(E)$. Then 
$\supp T\subseteq E$ and, therefore, $\langle T, F\rangle=0$ for any
$F\in\Phi(E)$. By Theorem~\ref{th6}, $T\in{\mathfrak M}_{min}(E)$, which gives
us the necessary inclusion 
${\mathfrak M}_{max}(E)\subseteq {\mathfrak M}_{min}(E)$.

$(ii)\Rightarrow (iv)$: if ${\mathfrak M}_{min}(E)={\mathfrak M}_{max}(E)$ then
$\text{bil }1\otimes{\mathfrak M}_{min}(E)=\text{bil }
1\otimes{\mathfrak M}_{max}(E)$ and
by Corollary~\ref{cor1} and Theorem~\ref{th7} we obtain that any
$E$-pair can be s-approximated  by simple $E$-pairs.

$(iv)\Leftrightarrow (v)$. We prove that the approximation  of operator-valued 
functions
in the strong operator topology in $B(l_2\otimes L_2(X,\mu))$ is equivalent 
to the
approximation almost everywhere in the strong operator topology in $B(l_2)$. 
In fact, let $P_n(x)$, $P(x)\in B(l_2\otimes L_2(X,\mu))$, $P_n(x)\to P(x)$ 
 almost everywhere on $(X,\mu)$ in the strong operator topology in
$B(l_2)$ and take $\varphi=\sum_{k=1}^N
\varepsilon_k(x)\vec{\xi_k}$, where $\varepsilon_k(\cdot)$ is the 
characteristic
function of a set of finite measure and $\vec{\xi_k}\in l_2$.  It easily follows from
 the Lebesgue theorem that $||P_n\varphi-P\varphi||\to 0$ as $n\to\infty$.
 Since the measure
$\mu$ is sigma-finite, the set of all such $\varphi$  is
dense in $l_2\otimes L_2(X,\mu)$. Therefore $||P_n\varphi-P\varphi||\to 0$, 
$n\to\infty$,
for any $\varphi\in  l_2\otimes L_2(X,\mu)$. 

If now a sequence, $\{P_n\}$,
of projection-valued functions converges to $P$ in the strong operator topology
in $B(l_2\otimes L_2(X,\mu))$, then there exists a subsequence converging
almost everywhere on $(X,\mu)$ in the strong operator topology in $B(l_2)$.
To see this  choose a dense set of vectors, $\{\vec{\xi_n}\}$, in $l_2$.
Then $$\int_{A}||P_n(x)\vec{\xi_k}-P(x)\vec{\xi_k}||d\mu(x)\to 0, \  
n\to\infty$$
for each $k$ and each measurable set $A$ of finite measure. 
Let $A_1\subseteq A_2\subseteq\ldots\subseteq A_n\subseteq\ldots $ be a sequence 
of sets
of finite measure such that $X=\cup_{j=1}^{\infty}A_j$.
By the Riesz theorem there exists a subsequence 
$\{P_{k1}\}_{k=1}^{\infty}$ such
that $\lim_{k\to\infty}P_{k1}(x)\vec{\xi_1}= P(x)\vec{\xi_1}$ a.e.\ on $A_1$. Then choose a subsequence $\{P_{k2}\}_{k=1}^{\infty}$
of $\{P_{k1}\}_{k=1}^{\infty}$ such that 
$\lim_{k\to\infty}P_{k2}(x)\vec{\xi_1}= P(x)\vec{\xi_1}$ a.e.\
on $A_2$. Proceeding in this fashion we obtain
a series of sequences
$$\{P_n\}_{n=1}^{\infty}\supset\{P_{k1}\}_{k=1}^{\infty}\supset
\{P_{k2}\}_{k=1}^{\infty}\supset\ldots\supset
\{P_{kj}\}_{k=1}^{\infty}\supset\ldots,$$
such that $\lim_{k\to\infty}P_{kj}(x)\vec{\xi_1}= P(x)\vec{\xi_1}$
almost everywhere on $A_j$. 

Consider now the diagonal sequence $\{P_{kk}\}_{k=1}^{\infty}$. 
Clearly $\lim_{k\to\infty}P_{kk}(x)\vec{\xi_1}= P(x)\vec{\xi_1}$ a.e.\ on each
$A_j$ and therefore on $X$. Set $P^{l1}=P_{ll}$, $l=1,2,\ldots$.
Using the same arguments we can
find a subsequence, $\{P^{l2}\}_{l=1}^{\infty}$, of 
$\{P^{l1}\}_{l=1}^{\infty}$ such that $\lim_{l\to\infty}P^{l2}(x)\vec{\xi_2}= 
P(x)\vec{\xi_2}$ a.e.\ on $X$ and then  $\{P^{lk}\}_{l=1}^{\infty}$, of 
$\{P^{l1}\}_{l=1}^{\infty}$ such that $\lim_{l\to\infty}P^{lk}(x)\vec{\xi_m}= 
P(x)\vec{\xi_m}$ a.e.\ on $X$ for any $m\leq k$ so that
$\lim_{l\to\infty}P^{ll}(x)\vec{\xi_k}= 
P(x)\vec{\xi_k}$ a.e.\ on $X$ for any $k$. Since $\{\vec{\xi_k}\}$ is dense 
in $l_2$ and the sequence $\{P^{ll}\}_{l=1}^{\infty}$ is bounded, 
$$\lim_{l\to\infty}P^{ll}(x)\vec{\xi}= 
P(x)\vec{\xi}  \text{ a.e. on  }X \text{ for any }\vec{\xi}\in l_2.$$

$(iv)\Rightarrow (ii)$: if $T\in{\mathfrak M}_{max}(E)$, we have
$\text{bil }1\otimes T\supseteq\overline{\{(P,Q):(P,Q)\text{ is a simple $E$ 
pair }\}}^s$, due to Theorem~\ref{th7}; $(iv)$ implies now 
$\text{bil }1\otimes T\supseteq \tilde S_E$ and hence 
$T\in{\mathfrak M}_{min}(E)$.
\end{proof}

\begin{remark}\rm
The equivalence $(i)\Leftrightarrow (iii)$ was essentially proved in 
\cite{arv} and $(i)\Leftrightarrow (ii)$ in \cite{davidson} but
using some other methods.
\end{remark}

We use the equivalence $(i)\Leftrightarrow (v)$ to obtain the 
following result.
\begin{theorem1}[{\bf Inverse Image Theorem}]\label{preimage}
Let $(X,\mu)$, $(Y,\nu)$, $(X_1,\mu_1)$ and $(Y_1,\nu_1)$ be standard
Borel spaces
with measures, $\varphi:X\mapsto X_1$, $\psi:Y\mapsto Y_1$ Borel mappings.
Suppose that the measures $\varphi_{*}\mu$, $\psi_{*}\nu$
are absolutely continuous with respect to the measures $\mu_1$ and $\nu_1$ 
respectively. If
a Borel set $E_1\subseteq X_1\times Y_1$ is a set of $\mu_1\times\nu_1$-synthesis
then $(\varphi\times\psi)^{-1}(E_1)$ is a set of $\mu\times\nu$ synthesis.
\end{theorem1}
\begin{proof}
To prove the theorem we will need to prove first an auxiliary lemma.

\begin{lemma1}\label{sec}
Let $(X,\mu)$, $(Y,\nu)$ be standard Borel spaces with measures and 
$f:X\to Y$ be a Borel map. 
Then there exists a $\nu$-measurable set $N\subset f(X)$, $\nu(N)=0$, 
such that $f(X)\setminus N$ is Borel and
if $u:X\to{\mathbb R}$ is a bounded Borel function then for any 
$\varepsilon>0$ there exists a Borel map $g:f(X)\setminus N\to X$ such that 
$f(g(y))=y$ for every $y\in f(X)\setminus N$ and 
$u(g(f(x)))>u(x)-\varepsilon$ a.e. on $X$.
\end{lemma1}

\begin{proof}
Assume first that the map $f:X\to Y$ is surjective. 
For any such map there exists a 
Borel section, i.e., a map $g:Y\to X$ which satisfies $f(g(y))=y$, $y\in Y$ 
(see, for example, \cite{takesaki}).
Since $u:X\to{\mathbb R}$ is bounded, $u(X)\subseteq[a,b]$. 
Let $a=a_0<a_1<\ldots <a_n=b$ be
a partition of $[a,b]$ such that $a_{i+1}-a_i<\varepsilon$. Set 
$$X_j=u^{-1}([a_j,a_{j+1})),\
Y_j=f(X_j),\ Y_j'=Y_j\setminus(\cup_{k>j} Y_k).$$ Then each $Y_j'$ is 
the image of $X_j'=X_j\setminus(\cup_{k>j}f^{-1}(Y_k))$. We have also that
$\cup_j Y_j'=Y$, $Y_i'\cap Y_j'=\emptyset$, $i\ne j$, and since 
every $Y_i'$ is an analytic space, we obtain that $Y_i'$ must be Borel (see, 
for example, \cite[Theorem~A.3]{takesaki}). 
Let $g_j:Y_j'\to X_j'$ be a Borel section for $f|_{X_j'}$. Then 
the functions $g_j$ determine a Borel section, $g$, for $f$.
Clearly, $g(Y_j)\subseteq\cup_{i\geq j}X_i$ so that 
$u(g(y))\geq a_j$ for each $y\in Y_j$
and therefore $u(g(f(x)))\geq a_j$ for any $x\in X_j$.
As $u(x)\in[a_j,a_{j+1})$ for $x\in X_j$, we obtain $u(g(f(x)))> 
u(x)-\varepsilon$ for each
$x_j\in X_j$ and therefore for each $x\in X$.

For the general case  consider the image $f(X)$ which is an analytic subset 
of $Y$. By \cite[Theorem~A.13]{takesaki} there exists a $\nu$-measurable set 
$N\subset f(X)$ of zero measure
such that $f(X)\setminus N$ is Borel. Set $\tilde X=f^{-1}(f(X)\setminus N)$. 
Then $f$ is a Borel map from the Borel set $\tilde X$ onto $f(X)\setminus N$. 
Thus, given  $\varepsilon>0$, there
exists a Borel map $g:f(X)\setminus N\to X$ such that $f(g(y))=y$ for every 
$y\in f(X)\setminus N$ and 
$u(g(f(x)))>u(x)-\varepsilon$ on $\tilde X$. Since 
$X\setminus\tilde X\subseteq f^{-1}(N)$, we have
that $\mu(X\setminus\tilde X)=0$ and the inequality holds almost everywhere 
on $X$.
\end{proof}

Set $E=(\varphi\times\psi)^{-1}(E_1)$. By Theorem~\ref{main}, we shall
have established the theorem if we prove  that any $E$-pair can be 
approximated a.e.\ in the strong operator topology of $B(l_2)$ by simple 
$E$-pairs.  Since, by Theorem~\ref{th7}, the approximated pairs form a 
bilattice it would be enough to prove that any
$E$-pair is majorized by an approximated pair. 

Let $(P,Q)$ be an $E$-pair. Choose a dense sequence $\xi_n$ in $l_2$ and a 
sequence $\e_n>0$, $\e_n\to 0$. 
Set $u_n(x)=(P(x)\xi_n,\xi_n)$. By Lemma~\ref{sec}, there are null sets 
$N_n\subset X_1$, $M_n\subset X$ and  a Borel map 
$g_n:\varphi(X)\setminus N_n\to X$, such that 
$\varphi(g_n(x_1))=x_1$, for $x_1\in \varphi(X)\setminus N_n$, and 
$u_n(g_n(\varphi(x)))>u_n(x)-\varepsilon_n$, for $x\in X\setminus M_n$.

For $x_1\in\varphi(X)\setminus N$, where $N=\cup_{n=1}^{\infty}N_n$,
set $$\hat P(x_1)=\bigvee_n P(g_n(x_1)).$$
Then for any $x\in X\setminus M$, where $M=\cup_{n=1}^{\infty} M_n$, one has
\begin{gather*}
(P(x)\xi_n,\xi_n)=u_n(x)<u_n(g_n(\varphi(x)))+\varepsilon_n=\\=
(P(g_n(\varphi(x)))\xi_n,\xi_n))+\e_n\leq
(\hat P(\varphi(x))\xi_n,\xi_n)
+\varepsilon_n.
\end{gather*}
It easily follows  that 
\begin{equation}\label{f1}
P(x)\leq\hat P(\varphi(x)),\quad x\in X\setminus M.
\end{equation}
Similarly, we construct null sets $M'\subset Y$, $N'\subset Y_1$, functions 
$g_n':\psi(Y)\setminus N'\to Y$ and set
$\hat Q(y_1)=\bigvee_n Q(g_n'(y_1))$ so that  
\begin{equation}\label{f2}
Q(y)\leq \hat Q(\psi(y)), \quad y\in Y\setminus M'.
\end{equation}
Thus $(P,Q)$ is majorized by $(\hat P\circ\varphi,\hat Q\circ\psi)$.

Setting $\hat P=0$ and $\hat Q=0$ on the complements of
$\varphi(X)\setminus N$ and $\psi(Y)\setminus N'$ respectively, we have
that $(\hat P,\hat Q)$ is an $E_1$-pair. Indeed, let $(x_1,y_1)\in E_1$,
$x_1\in\varphi(X)\setminus N$, $y_1\in \psi(Y)\setminus N'$, then 
$$P(g_n(x_1))\perp Q(g_m'(y_1))$$
for any $n$, $m$. Hence
$$\hat P(x_1)\perp\hat Q(y_1).$$
It follows that there are simple $E_1$-pairs $(\hat P_n,\hat Q_n)$ with
$\hat P_n(x_1)\to\hat P(x_1)$ a.e. ($x_1\notin S$), 
$\hat Q_n(y_1)\to\hat Q(y_1)$ a.e. ($y_1\notin S'$). Let
$$P_n(x)=\hat P_n(\varphi(x)),\quad Q_n(y)=\hat Q_n(\psi(y)).$$
Then $P_n(x)\to\hat P(\varphi(x))$ a.e., $Q_n(y)\to\hat Q(\psi(y))$ a.e.
Indeed, let $\tau=\{x:\varphi(x)\in S\}$, then
$$\mu(\tau)=\mu(\{x:\varphi(x)\in S\})=\varphi_*\mu(S)=0,$$
because $\varphi_*\mu$ is absolutely continuous with respect to $\mu_1$. Similarly,
$\nu(\tau')=0$, where $\tau'=\{y:\psi(y)\in S'\}$. This shows that the pair
$(\hat P\circ\varphi,\hat Q\circ\psi)$ is approximable by simple pairs.
The proof is complete.
\end{proof}

\begin{cor1}\label{measures}
Let $E\subseteq X\times Y$ be a set of synthesis with respect to a pair of
measures $(\mu_1,\nu_1)$, $\mu_1\in M(X)$, $\nu_1\in M(Y)$. Then $E$ is a set 
of
$(\mu,\nu)$-synthesis for any $\mu\in M(X)$, $\nu\in M(Y)$ such that
$\mu\leq \mu_1$, $\nu\leq \nu_1$. 
\end{cor1}

\begin{proof}
Follows from Theorem~\ref{preimage} applied to the identity mappings 
$\varphi$ and $\psi$.
\end{proof}

 Suppose that $f_i$ and $g_i$, 
$i=1,\ldots,n$, are Borel maps of   standard Borel spaces $(X,\mu)$ and 
$(Y,\nu)$ into an ordered standard Borel space $(Z,\leq)$. Then the set 
$E=\{(x,y)\mid f_i(x)\leq g_i(y), 
i=1,\ldots,n\}$ is called a set of width $n$.

\begin{theorem1}\label{fw}
Any set of finite width is synthetic with respect to the measures
$\mu$, $\nu$.
\end{theorem1}
\begin{proof}
Let $E$ be a set of width $n$, i.e.
$E=\{(x,y)\in X\times Y\mid f_i(x)\leq g_i(y), i=1,\ldots,n\}$, where
$f_i:X\to Z$, $g_i:Y\to Z$ are Borel functions.
We define  mappings $F:X\to Z^n$ and $G:Y\to Z^n$ by setting
$F(x)=(f_1(x),\ldots,f_n(x))$, $G(y)=(g_1(y),\ldots,g_n(y))$.
Put $\mu_1=F_{*}\mu$, $\nu_1=G_{*}\nu$. 
Let $E_1=\{(x,y)\in Z^n\times Z^n\mid x_i\leq y_i,i=1,\ldots,n\}$.
By \cite{arv}, 
$E_1$ is a set of 
$\mu_1\times\nu_1$-synthesis if the measures $\mu_1$ and $\nu_1$ are equal. 
In general, consider the measure $\lambda=\mu_1+\nu_1$, then we can conclude
that $E_1$ is a set
of $\lambda\times\lambda$-synthesis and applying now 
Corollary~\ref{measures} we obtain that $E_1$ is a set of synthesis
with respect to $\mu_1$, $\nu_1$.
 It follows now from Theorem~\ref{preimage} that
$(F\times G)^{-1}(E_1)=E$ is a set of $\mu\times\nu$-synthesis.
\end{proof}

\begin{remark}\rm
Arveson, \cite{arv}, introduced the class of finite width lattices as 
those which are generated by a finite set of nests (linearly ordered
lattices). He proved that all finite width lattices are synthetic.
Todorov, \cite{todorov}, defined a subspace map (see \cite{erdos}) of
finite width and proved that such subspace maps are synthetic. This
result is in fact equivalent to our, actually a subspace map is a counterpart 
of a bilattice.
Synthesizability of special sets of width two
(``nontriangular'' sets) was proved in \cite{kt, sh}.
\end{remark}

In \cite{arv}[Problem, p.487] Arveson also posed a question whether or not 
the lattice generated by a synthetic lattice and a lattice of finite width 
is synthetic. 
Next result shows that the answer is no. The example we construct is 
inspired by the Varopoulos example (\cite{var1}) of a set of spectral 
synthesis for the Fourier algebra $A({\mathbb R}^2)$ whose intersection 
with a subgroup does not admit synthesis.

Let $F$ denote the Fourier transform in ${\mathbb R}^n$ and let 
$A({\mathbb R}^n)$ be the Fourier algebra $FL_1({\mathbb R}^n)$ which is
a Banach algebra with the norm $||Ff||_{A}=||f||_{L_1}$. Recall
that a closed set $K\subseteq{\mathbb R}^n$ admits spectral synthesis for
$A({\mathbb R}^n)$
if for every $f\in A({\mathbb R}^n)$ vanishing
on $K$ there exists a sequence $f_n\in A({\mathbb R}^n)$ such that 
$f_n$ vanishes on an open set containing $K$ and $||f_n-f||_A\to 0$
as $n\to\infty$.

A commutative lattice ${\mathcal L}$ is called synthetic if the only
ultra-weakly closed algebra ${\mathcal A}$ satisfying 
$\text{lat }{\mathcal A}={\mathcal L}$ and ${\mathcal L}'\subseteq{\mathcal A}$
is the algebra $\text{alg }{\mathcal L}$. If ${\mathcal L}_1$,
${\mathcal L}_2$ are two lattices we will denote by ${\mathcal L}_1\bigvee
{\mathcal L}_2$ the lattice generated by ${\mathcal L}_1$ and ${\mathcal L}_2$.
\begin{theorem1}\label{synfw}
There exist a synthetic lattice ${\mathcal L}_1$ and a lattice ${\mathcal L}_2$
of finite
width such that ${\mathcal L}_1\bigvee{\mathcal L}_2$ is not synthetic.
\end{theorem1}
\begin{proof}
Let $G\subset{\mathbb R}$ be a set which does not admit spectral synthesis
for $A({\mathbb R})$. Set
$$E=\{(x,t):d(x,G)\leq t\}\subset{\mathbb R}^2.$$
Here $d(x,G)$ denotes the distance between $x$ and $G$. Then $E$ is
a set of spectral synthesis for $A({\mathbb R}^2)$. Indeed, if 
$f(x,t)\in A({\mathbb R}^2)$ vanishes on $E$ then $f_n(x,t)=f(x,t+1/n)\in
A({\mathbb R}^2)$ vanishes on $E_n=\{(x,t)\mid d(x,G)< t+1/n\}$ containing $E$
and $||f_n-f||_A\to 0$ as $n\to\infty$.

The intersection $E\cap({\mathbb R}\times\{0\})=G\times\{0\}$ does
not admit spectral synthesis. In fact, otherwise, given 
$f(x,t)\in A({\mathbb R}^2)$, $f(x,0)=0$ for $x\in G$, there exists
a sequence $f_n(x,t)\in A({\mathbb R}^2)$ such that $f_n(x,t)=0$ on nbhd of 
$G\times\{0\}$ and $||f_n(x,t)-f(x,t)||_A\to 0$, $n\to\infty$. 
Now it is enough to
see that $f_n(x,0)$, $f(x,0)\in  A({\mathbb R})$, each $f_n(x,0)$ vanishes
on a nbhd of $G$ and $||f_n(x,0)-f(x,0)||_A\to\infty$ as $n\to\infty$, 
contradicting the assumption that $G$ is not a set of spectral synthesis.

If $m$ denotes the Lebesgue measure on ${\mathbb R}^2$, by \cite{F} 
we have that $$E^*=\{(x,y)\in{\mathbb R}^2\times{\mathbb R}^2\mid x-y\in E\}$$
is a set of $m\times m$-synthesis while
$$
(G\times\{0\})^*=\{(x,y)\in{\mathbb R}^2\times{\mathbb R}^2\mid x-y\in 
G\times\{0\}\}=E^*\cap L,
$$
where $L=\{(x,y)\in{\mathbb R}^2\times{\mathbb R}^2\mid x_2=y_2\}$,
is not $m\times m$-synthetic.

Let ${\mathcal L}$ and ${\mathcal L}_1$ be the lattices of projections
$P\oplus (1-Q)$, where $(P,Q)$ belongs to the bilattices 
$S_{(G\times\{0\})^*}$ and $S_{E^*}$ respectively.
Then ${\mathcal L}_1$ is synthetic while ${\mathcal L}$ is not.
In fact, if ${\mathcal A}$ is an ultra-wekly closed algebra such that
$\text{lat }{\mathcal A}={\mathcal L}$ 
and ${\mathcal L}'\subseteq{\mathcal A}$ ($\text{lat }{\mathcal A}
={\mathcal L}_1$ and ${\mathcal L}_1'\subseteq{\mathcal A}$) then 
${\mathcal A}=\left(\begin{array}{cc}{\mathcal A}_{11}&0\\
{\mathcal A}_{21}&{\mathcal A}_{22}\end{array}\right)$ where
${\mathcal A}_{ii}=L_{\infty}({\mathbb R}^2)$, ${\mathcal A}_{21}$
is an ultra-weakly closed $L_{\infty}({\mathbb R}^2)\times L_{\infty}
({\mathbb R}^2)$-bimodule such that $\text{bil }{\mathcal A}_{21}=
S_{(G\times\{0\})^*}$ ($\text{bil }{\mathcal A}_{21}=S_{E^*}$). 
The statement now follows from the synthesizability of $E^*$ and the 
non-synthesizability of $(G\times\{0\})^*$.

 Let $P_{\Sigma}$
denote the multiplication operator by the characteristic function
of the set $\Sigma$ and let ${\mathcal L}_2$ be the lattice of projections 
$P_{\Sigma}\oplus P_{\Sigma}$, 
where $\Sigma={\mathbb R}\times K$ and $K$ is a Borel subset of ${\mathbb R}$
($\Sigma$ is an increasing set for the partial ordering $x\leq y$, 
$x, y\in{\mathbb R}^2$ iff $x_2=y_2$). Then ${\mathcal L}_2$ is a set of 
width $2$ generated by the nests ${\mathcal C}$ and ${\mathcal C}^{\perp}$, 
where ${\mathcal C}=\{P_{\Sigma_t}\oplus P_{\Sigma_t}\mid
\Sigma_t={\mathbb R}\times[t,+\infty)\, t\in{\mathbb R}\}$.

What is left to prove is that 
${\mathcal L}={\mathcal L}_1\bigvee{\mathcal L}_2$.
Since $(G\times\{0\})^*=E^*\cap L$,
one easily sees that ${\mathcal L}_1$, ${\mathcal L}_2\subseteq{\mathcal L}$ and
therefore ${\mathcal L}_1\bigvee{\mathcal L}_2\subseteq{\mathcal L}$.
For the reverse inclusion we use the reflexivity of the CSL
${\mathcal L}_1\bigvee{\mathcal L}_2$.
 Direct  verification shows that
$$\text{alg }{\mathcal L}_1\bigvee{\mathcal L}_2=\{(T_{ij})_{i,j=1}^2
\mid T_{11}, T_{22}\in 
L_{\infty}({\mathbb R}^2), T_{12}=0, \text{supp }T_{21}\subseteq E^*\cap L\}.$$
Therefore if $P\oplus (1-Q)\in{\mathcal L}$, i.e. $P=P_{\alpha}$, $Q=P_{\beta}$
for some Borel sets $\alpha$, $\beta$ such that 
$(\alpha\times\beta)\cap(E^*\cap L)=\emptyset$,
we have 
$P\oplus(1-Q)\in\text{lat }\text{alg }{\mathcal L}_1\bigvee{\mathcal L}_2
={\mathcal L}_1\bigvee{\mathcal L}_2$ and ${\mathcal L}\subseteq 
{\mathcal L}_1\bigvee{\mathcal L}_2$.
 
\end{proof}
\section{General bilattices}\label{sssec122}

Let $h_0$ be a   function on $[0,1]$ defined by
$h_0(0)=0$ and $h_0(t)=1$ for
$t\ne 0$, and let $h_1(t)=1-h_0(1-t)$. It is clear that for any
positive contraction $A$, $h_0(A)$ is the projection onto the range of $A$,
$h_1(A)$ is the projection onto the subspace of invariant vectors.
It is easy to see (for example, approximating $h_1(t)$ by $t^{\alpha}$,
$\alpha\to 0$) that $h_i$ are operator monotone, i.e., if
$A$, $B\in B(H)$, $0\leq A\leq B\leq 1$, then $h_i(A)\leq h_i(B)$.

Recall, given a commutative ${\mathcal D}_1\times{\mathcal D}_2$-bilattice
$S$, 
$$F_S=\{(A,B)\in (B(l_2)\bar\otimes {\mathcal D}_1)\times
(B(l_2)\bar\otimes {\mathcal D}_2):(L_{\varphi}(A), L_{\varphi}(B))
\in\text{Conv }S\text{ for any }\varphi\}.$$

\begin{lemma1}\label{functh}
Let ${\mathcal D}_1$, ${\mathcal D}_2$ be  commutative von Neumann algebras 
on Hilbert spaces
$H_1$, $H_2$ and let $S$ be a  ${\mathcal D}_1\times {\mathcal D}_2$-bilattice.
Then, for any $(A,B)\in F_S$, 
$(h_0(A), h_1(B))\in F_S$ and $(h_1(A), h_0(B))\in F_S$.
\end{lemma1}
\begin{proof} If ${\mathcal D}_1$, ${\mathcal D}_2$ are masas in 
separable spaces $H_1$, $H_2$,
then the assertion  follows from Lemma~\ref{cor0a}. 
Indeed, if
$A(x)+B(y)\leq 1$, then $$h_0(A(x))\leq h_0(1-B(y))=1-h_1(B(y))$$
and $(h_0(A), h_1(B))\in F_S$. Similarly, $(h_1(A), h_0(A))\in F_S$.

Assume now that ${\mathcal D}_1$, ${\mathcal D}_2$ are arbitrary commutative
von Neumann algebras acting on separable Hilbert spaces. Let $x_1$ and $x_2$
be separating vectors for ${\mathcal D}_1$ and ${\mathcal D}_2$, and let 
$K_i=\overline{[{\mathcal D}_ix_i]}$, $i=1,2$. Then the restriction of
$B(l_2)\bar\otimes {\mathcal D}_i$ to $l_2\otimes K_i$ is injective. Now,
since the restriction of ${\mathcal D}_i$ to $K_i$ is a masa and the 
restriction
of $(A,B)\in F_S$  to $(l_2\otimes K_1)\times(l_2\otimes K_2)$ 
belongs to $F_{\bar S}$, where ${\bar S}$ is the restriction of $S$ to 
$K_1\times K_2$, the problem is reduced to the above.

Furthermore, the statement is true when ${\mathcal D}_1$, ${\mathcal D}_2$ 
are countably generated.
To see this it is enough to prove that if $x_1,\ldots,x_n$
and $y_1,\ldots,y_n$ are vectors in $l_2\otimes H_1$ and $l_2\otimes
H_2$, then there exist a pair $(C,D)\in F_S$
such that $h_0(A)x_i=Cx_i$ and $h_1(B)y_i=Dy_i$, $i=1,\ldots,n$.
If $x_k=(x_{kj})$, $y_k=(y_{kj})$, $x_{kj}\in H_1$, $y_{kj}\in H_2$,
 we define $K_1$ and $K_2$
to be the closed linear spans of vector $Xx_{kj}$, $X\in {\mathcal D}_1$, and
$Yy_{kj}$, $Y\in {\mathcal D}_2$, respectively. Then $K_1$ and $K_2$ are
separable and we come to the previous case.
 
Now, to prove the assertion in general situation, it is sufficient to
show that each ${\mathcal D}_i$ contains a countably generated von Neumann 
algebra,
$\hat{\mathcal D}_i$, such that $(A,B)\in F_{\hat S}$, where $\hat S$ is the 
intersection
of $S$ with $\hat {\mathcal D}_1\times \hat {\mathcal D}_2$. For this take  
a dense sequence of unit vectors, $\{\xi_n\}$, in $l_2$. 
For each pair $(L_{\xi_n}(A),L_{\xi_n}(B))$ there exists
a sequence, $(A_k^n, B_k^n)$, from the convex linear span, 
$\text{conv }S$, of $S$, which
converges to the pair uniformly. Let $S'$ be the set of all pairs
of projections $(p,q)\in S$ which participate in the linear combinations
for $(A_k^n, B_k^n)$. Then $\hat {\mathcal D}_1$ and $\hat {\mathcal D}_2$ 
can be defined  as von Neumann algebras generated by $\pi_1(S')$ and 
$\pi_2(S')$, $\pi_i$ being the projection onto the $i$-th coordinate.
\end{proof}
\begin{lemma1}
$\tilde S$ is a billatice.
\end{lemma1}
\begin{proof}
Let $(P,Q)\in\tilde S$ and 
$P_1\in{\mathcal P}_{B(l_2)\bar\otimes {\mathcal D}_1}$, 
$Q_1\in{\mathcal P}_{B(l_2)\bar\otimes {\mathcal D}_2}$,
$P_1\leq P$, $Q_1\leq Q$. Then $L_{\varphi}(P_1)\leq L_{\varphi}(P)$, 
$L_{\varphi}(Q_1)\leq L_{\varphi}(Q)$ for each state $\varphi$ on $B(l_2)$ 
so that
$$E_{L_{\varphi}(P_1)}([a,1])\leq E_{L_{\varphi}(P)}([a,1])\text{ and }
E_{L_{\varphi}(Q_1)}([b,1])\leq E_{L_{\varphi}(Q)}([b,1])$$
for any $0\leq a,b\leq 1$. Applying now Lemma~\ref{th3a2} we obtain 
$(P_1,Q_1)\in\tilde S$.

That $\tilde S$ is closed under the operations $(\bigvee,\bigwedge)$, 
$(\bigwedge,\bigvee)$ follows from
\begin{eqnarray*}
(P_1\bigvee P_2,Q_1\bigwedge Q_2)=(h_0((P_1+P_2)/2),h_1((Q_1+Q_2)/2),\\
(P_1\bigwedge P_2,Q_1\bigvee Q_2)=(h_1((P_1+P_2)/2),h_0((Q_1+Q_2)/2)
\end{eqnarray*} 
and the previous lemma.
\end{proof}

Our next goal is to show that $\tilde S$ is reflexive. We will deduce this 
from a general criteria of reflexivity. To formulate it we need some
definitions and notations.

Let $S$ be an ${\mathcal R}\times{\mathcal R}$-bilattice, where
${\mathcal R}$ is a von Neumann algebra on a Hilbert space $H$, and let 
${\mathcal M}$
be a von Neumann algebra on $H$. Denote by $I({\mathcal M})$ the semigroup 
of all isometries in ${\mathcal M}$. We say that $S$ is 
{\it ${\mathcal M}$-invariant} if
\begin{itemize}
\item
 $S$ contains all pairs $(P,1-P)$, $P\in{\mathcal P}_{\mathcal M}$.

\item
 If $U\in I(\M)$ then a pair $(P,Q)\in{\R}\times{\R}$ belongs to $S$
if and only if $(UPU^*,UQU^*)$ belongs to $S$.
\end{itemize}
For any bilattice $S$ we set
$$\Omega_S=\{(x,y)\in H\times H\mid \exists (P,Q)\in S\text{ with } 
Px=x, Qy=y\}.$$
If $S$ is clear we write $\Omega$ instead of $\Omega_S$.
A bilattice $S$ is called {\it stable} if $\Omega_S$ is norm-closed in 
$H\oplus H$.
\begin{theorem1}\label{threfl}
Any  ${\mathcal R}\times{\mathcal R}$-bilattice which is stable and invariant 
with respect to a properly
infinite von Neumann algebra is reflexive.
\end{theorem1}
\begin{proof}
Suppose that $S$ is stable and $\M$-invariant, where $\M$ is properly infinite.
Note first that ${\mathfrak M}(S)\subseteq\M'$. Indeed, if $T\in{\mathfrak M}(S)$
then $(1-P)TP=0$ for any $P\in{\mathcal P}_{\M}$, and similarly $PT(1-P)$, 
hence $TP=PT$ and $T\in\M'$, because ${\mathcal P}_{\M}$ generates $\M$.

\vspace{0.2cm}

{\bf Claim 1.} Let $U\in I(\M)$. If $(U^*x,y)\in\Omega$ then $(x,Uy)\in\Omega$.

Indeed, let $(P,Q)\in S$ such that $PU^*x=U^*x$, $Qy=y$. Consider
$P_1=UPU^*$, $Q_1=UQU^*$. Then $P_1x=UU^*x$, $Q_1Uy=Uy$ and thus 
$Uy\in Q_1H\cap UU^*H$. Set 
$$P_2=P_1\bigvee (1-UU^*),\quad Q_2=Q_1\bigwedge UU^*.$$
Then $P_2H$ contains $UU^*x$ and $(1-UU^*)x$, hence $P_2H$ contains $x$,
i.e. $P_2x=x$. On the other hand $Q_2H$ contains $Uy$. So $Q_2Uy=Uy$.
Clearly, $(P_2,Q_2)\in S$ and we get $(x,Uy)\in\Omega$. 

\vspace{0.1cm}

Now we prove the converse statement.

\vspace{0.2cm}

{\bf Claim 2.} If $(x,Uy)\in\Omega$, $U\in I(\M)$ then $(U^*x,y)\in \Omega$.

Indeed, let $(P,Q)\in S$, $Px=x$, $QUy=Uy$. Set $$P_1=P\bigvee (1-UU^*),
Q_1=Q\bigwedge UU^*.$$ Then $P_1x=x$, $Q_1Uy=Uy$,
$P_1\geq 1-UU^*$, $Q_1\leq UU^*$. It follows that $P_1$, $Q_1$ commute 
with $UU^*$. Hence $P_2=U^*P_1U$ and $Q_2=U^*Q_1U$ are projections.
To see that $(P_2,Q_2)\in S$ note that 
$(UP_2U^*,UQ_2U^*)=(UU^*P_1,UU^*Q_1)\in S$, since $UU^*P_1\leq P_1$, 
$UU^*Q_1\leq Q_1$.
It remains to show that $P_2U^*x=U^*x$ and
$Q_2y=y$. Indeed, 
\begin{eqnarray*}
&P_2U^*x=U^*P_1UU^*x=U^*UU^*P_1x=U^*UU^*x=U^*x,\\
&Q_2y=U^*Q_1Uy=U^*Uy=y.
\end{eqnarray*}

Our claim is proved.

\vspace{0.2cm}
For $(x,y)\in H\times H$, we denote by $v_{x,y}$ the restriction of
the vector state $w_{x,y}$ to $\M'$.

{\bf Claim 3.} If $(x,y)\in\Omega$, $v_{x,y}=v_{x,z}$ then $(x,z)\in\Omega$.

To show this set $t=y-z$. Then $v_{x,t}=0$, $\M'x\perp\M't$. Defining
$R$ to be the projection onto $\overline{\M'x}$ we have $R\in\M$, $Rx=x$
and $(1-R)t=t$.

Let now $(P,Q)\in S$, $Px=x$, $Qy=y$. Set $P_1=P\bigwedge R$, 
$Q_1=Q\bigvee (1-R)$. Then $(P_1,Q_1)\in S$, $P_1x=x$, $Q_1z=Q_1(y-t)=y-t=z$.
We proved that $(x,z)\in\Omega$.

\vspace{0.1cm}

Since $\M$ is properly infinite there are $U_1$, $U_2\in I(\M)$ with
$U_1H\perp U_2H$. We fix such a pair of isometries.

\vspace{0.2cm}

{\bf Claim 4.} If $(x_1,y_1)\in\Omega$ and $v_{x_1,y_1}=v_{x_2,y_2}$ then 
$(x_2,y_2)\in\Omega$.

Indeed, set $x=U_1x_1+U_2x_2$. Then $x_1=U_1^*x$. 
Hence $(U_1^*x,y_1)\in\Omega$.
By Claim 1, $(x,U_1y_1)\in\Omega$. Since 
$$v_{x,U_1y_1}=v_{U_1^*x,y_1}=v_{x_1,y_1}=v_{x_2,y_2}=v_{x,U_2y_2},$$
we obtain from Claim 3 that $(x,U_2y_2)\in\Omega$. Now by Claim 2,
$(U_2^*x,y_2)\in\Omega$, that is $(x_2,y_2)\in\Omega$. The claim is proved.

\vspace{0.1cm}

Set now $$W=\{v_{x,y}\mid (x,y)\in\Omega\}.$$

\vspace{0.2cm}

{\bf Claim 5.} $W$ is a linear subspace in the space $(\M')_*$ of all
$\sigma$-weakly continuous functionals on $\M'$.

Indeed,
$$v_{x_1,y_1}+v_{x_2,y_2}=v_{x,U_1^*y_1}+v_{x,U_2^*y_2}=v_{x,y},$$
where $x=U_1x_1+U_2x_2$. We know from the preceding claim
that $(x,U_1^*y_1)$ and $(x,U_2^*y_2)$ belong to $\Omega$. Let $(P_1, Q_1)
\in S$,  $(P_2,Q_2)\in S$ such that $$P_1x=x,\  Q_1U_1^*y_1=U_1^*y_1,\ 
P_2x=x,\  Q_2U_2^*y_2=U_2^*y_2.$$ Then setting
$P=P_1\bigwedge P_2$, $Q=Q_1\bigvee Q_2$ we have $Px=x$, $Qy=y$. Thus
$(x,y)\in\Omega$ and $W+W\subseteq W$.

\vspace{0.2cm}

{\bf Claim 6.} $W$ is norm-closed.

Let $\varphi_n\to\varphi$, $\varphi_n\in W$. Since $\varphi$ is 
$\sigma$-weakly continuous and $\M'$ has a separating vector, $\varphi=v_{x,y}$
for some $x\in H$, $y\in H$. Since $\M'$ has the properly infinite commutant,
there are $x_n$, $y_n\in H$ such that $\varphi_n=v_{x_n,y_n}$, 
$||x_n-x||\to 0$, $||y_n-y||\to 0$ (\cite{Sh}). By Claim 4,
$(x_n,y_n)\in\Omega$. Since $S$ is stable, $(x,y)\in\Omega$ and
$\varphi\in W$. We proved that $W$ is norm-closed.

\vspace{0.1cm}

Recall that $\M'$ is the dual of $(\M')_*$. So for ${\mathcal A}\subseteq\M'$, 
${\mathcal B}\subseteq(\M')_*$ we write
$${\mathcal A}_{\perp}=\{\varphi\in(\M')_*\mid {\mathcal A}
\subseteq\ker\varphi\}, \quad
{\mathcal B}^{\perp}=\{T\in\M\mid \varphi(T)=0, \forall 
\varphi\in{\mathcal B}\}.
$$
By the usual duality argument, $({\mathcal B}^{\perp})_{\perp}$ coincides
with the norm closure of ${\mathcal B}$, for any linear subspace 
${\mathcal B}\subseteq(\M')_*$.

\vspace{0.2cm}

{\bf Claim 7.} $W=({\mathfrak M}(S))_{\perp}$.

Indeed, suppose that $T\in W^{\perp}$. Then for any $(P, Q)\in S$,
$QTP=0$, because 
$$(QTPx,y)=w_{Px,Qy}(T)=0.$$
Thus $W^{\perp}={\mathfrak M}(S)$ and, by duality, 
$W={\mathfrak M}(S)_{\perp}$,
since $W$ is closed.

\vspace{0.2cm}

Now we can finish the proof of the theorem.

If $(P_0,Q_0)\in\text{bil }{\mathfrak M}(S)$ then $w_{P_0x,Q_0y}(T)=0$
for any $T\in{\mathfrak M}(S)$. Hence
$w_{P_0x,Q_0y}\in{\mathfrak M}(S)_{\perp}=W$. On the other hand,
for any $x\in P_0H$, $y\in Q_0H$ there are $(P_{x,y}, Q_{x,y})\in S$ with
$x\in P_{x,y}H$, $y\in Q_{x,y}H$.
Set
$$P_x=\bigwedge_{y\in Q_0H}P_{x,y},\quad Q_x=\bigvee_{y\in Q_0H} Q_{x,y}.$$
Then $x\in P_xH$, $Q_0H\subseteq Q_xH$. Let
$$P=\bigvee _{x\in P_0H} P_x, \ Q=\bigwedge_{x\in P_0H} Q_x,$$ then
$(P,Q)\in S$, $P_0\leq P$, $Q_0\leq Q$ implying $(P_0,Q_0)\in S$.
\end{proof}

Let $S$ be an ${\mathcal R}_1\times{\mathcal R}_2$-bilattice and
let $B_S$ denote the $({\mathcal R}_1\oplus{\mathcal R}_2)
\times ({\mathcal R}_1\oplus{\mathcal R}_2)$-bilattice  
generated by all pairs $(P\oplus(1-Q), (1-P)\oplus Q)$, where $(P,Q)\in S$.
It is easy to see that $B_S$ consists of all pairs $(P_1\oplus P_2,
Q_1\oplus Q_2)$, where $(P_1,Q_2)\in S$ and $Q_1\leq 1-P_1$, $P_2\leq 1-Q_2$.

\begin{proposition1}\label{prrefl}
An ${\mathcal R}_1\times{\mathcal R}_2$-bilattice $S$ 
is reflexive if and 
only if the bilattice $B_S$  is reflexive.
\end{proposition1}
\begin{proof}
Since $S$ is a bilattice, $(P,0)$, $(0,Q)\in S$ for any $P\in{\mathcal R}_1$,
$Q\in{\mathcal R}_2$. This implies 
\begin{eqnarray*}
\begin{array}{lll}
{\mathfrak M}(B_S)&=&\{(T_{ij})_{i,j=1}^2\mid (1-P)T_{11}P=QT_{22}(1-Q)=
QT_{21}P=0,\\
&&\\
&&(1-P)T_{12}(1-Q)=0, \forall (P,Q)\in S\}=\\
&&\\
&=&\{(T_{ij})_{i,j=1}^2\mid 
T_{ii}\in{\mathcal R}_i', i=1,2, T_{21}\in {\mathfrak M}(S), T_{12}=0\}
\end{array}
\end{eqnarray*}
and
\begin{eqnarray*}
\text{bil }{\mathfrak M}(B_S)=\{ (P_1\oplus P_2,Q_1\oplus Q_2)\mid
Q_iT_{ii}P_i=Q_2T_{21}P_1=0, \forall T=(T_{ij})_{i,j=1}^2\in
{\mathfrak M}(B_S)\}\\=\{ (P_1\oplus P_2,Q_1\oplus Q_2)\mid Q_1P_1=Q_2P_2=0,
(P_1,Q_2)\in \text{bil }{\mathfrak M}(S)\}
\end{eqnarray*}
giving the statement.
\end{proof}

Let now $S$ be again a commutative bilattice in
${\mathcal D}_1\times{\mathcal D}_2$ and let $\tilde S$ be the bilattice 
defined above.

\begin{theorem1}\label{tildeS}
The bilattice $\tilde S$ is reflexive.
\end{theorem1}
\begin{proof}
By Proposition~\ref{prrefl} and Theorem~\ref{threfl} it is sufficient to prove
that the bilattice $B_{\tilde S}$ is stable and $B(l_2)\otimes 1$-invariant.

Let $(x_n^1\oplus x_n^2, y_n^1\oplus y_n^2)\in\Omega_{B_{\tilde S}}$,
$x_n^i$, $y_n^i\in l_2\otimes H_i$, $i=1,2$,
 and $x_n^i\to x_i$, $y_n^i\to y_i$
as $n\to\infty$.
Then $p_n^ix_n^i=x_n^i$, $q_n^iy_n^i=y_n^i$ for some
$(p_n^1\oplus p_n^2,q_n^1\oplus q_n^2)\in B_{\tilde S}$.
We have $(p_n^1,q_n^2)\in\tilde S$ and $q_n^1\leq 1-p_n^1$,
$p_n^2\leq 1-q_n^2$.
We can also assume that the sequences $\{p_n^i\}$, $\{q_n^i\}$ are weakly 
convergent:
$$p_n^i\to a_i,\quad q_n^i\to b_i.$$
Clearly, $a_ix_i=x_i$, $b_iy_i=y_i$ and $b_1\leq 1-a_1$, $a_2\leq 1-b_2$.
Let  $P_i=h_1(a_i)$ and $Q_i=h_1(b_i)$ be the projections onto invariant
 vectors of
$a_i$ and $b_i$, $i=1,2$. It is easy to check that $(a_1,b_2)\in F_{\tilde S}$.
By Lemma~\ref{th3a2}, $(P_1,Q_2)\in\tilde S$. Moreover,
$Q_1\leq 1-P_1$, $P_2\leq 1-Q_2$. Thus
$(P_1\oplus P_2,Q_1\oplus Q_2)\in B_{\tilde S}$, 
$(x_1\oplus x_2,y_1\oplus y_2)\in \Omega_{B_{\tilde S}}$ and $B_{\tilde S}$ is
stable.

In order to prove $B(l_2)\otimes 1$-invariance we note first that
for any unit vector $\xi\in l_2$, any $u\in I(B(l_2))$ and 
$P\in B(l_2)\bar\otimes {\mathcal D}_i$, $i=1,2$,
$$L_{\xi}(P)=L_{u\xi}((u\otimes 1)P(u\otimes 1)^*)\text{ and }
L_{\xi}((u\otimes 1)P(u\otimes 1)^*)=L_{u^*\xi}(P)$$
implying that 

\begin{equation}\label{iff}
(P,Q)\in\tilde S \text{ iff }
((u\otimes 1)P(u\otimes 1)^*, (u\otimes 1)Q(u\otimes 1)^*)\in\tilde S.
\end{equation}
Since $u$ is an isometry, we have also that for any 
$P\in B(l_2)\bar\otimes {\mathcal D}_i$
$$P\leq 1-Q\Leftrightarrow (u\otimes 1)P(u\otimes 1)^*\leq 1-(u\otimes 1)Q
(u\otimes 1)^*$$
From this and (\ref{iff}) it follows that
$(P_1\oplus P_2,Q_1\oplus Q_2)\in B_{\tilde S}$ if and only if
$((u\otimes 1)(P_1\oplus P_2)(u\otimes 1)^*, 
(u\otimes 1)(Q_1\oplus Q_2)(u\otimes 1)^*)\in B_{\tilde S}$. 

Since for a state $\varphi$ on $B(l_2)$  and $p\in{\mathcal P}_{B(l_2)}$,
$$(L_{\varphi}(p\otimes 1), L_{\varphi}((1-p)\otimes 1)=
(\varphi(p),1-\varphi(p))=\varphi(p)(1,0)+(1-\varphi(p))(0,1)\in 
\text{Conv }S,$$
we have also 
$$(p\otimes 1\oplus p\otimes 1,(1-p)\otimes 1\oplus 
(1-p)\otimes 1)\in B_{\tilde S}.$$
We proved therefore that $B_{\tilde S}$ is $B(l_2)\otimes 1$-invariant.
\end{proof}

{\bf Proof of Theorem~\ref{smallest}.}
Since $\text{bil }{\mathfrak M}_0(S)\supseteq S$, we have only to prove 
the reverse inclusion.
Let $(P,Q)\in \text{bil }{\mathfrak M}_0(S)$. Then
$(1\otimes P, 1\otimes Q)\in \text{bil }{\mathfrak M}(\tilde S)$.
By Theorem~\ref{tildeS}, $(1\otimes P, 1\otimes Q)\in\tilde S$ and therefore
$(P,Q)\in S$. \qed

\section{ Operator synthesis and spectral synthesis} 

We recall first the definition of a set of spectral synthesis.
Let ${\mathcal A}$ be a unital semi-simple regular commutative Banach algebra 
with spectrum $X$, which is thus a compact Hausdorff space.  
We will identify ${\mathcal A}$ with a subalgebra of the algebra $C(X)$ of 
continuous 
complex-valued functions on $X$ in our notation.
If $E\subseteq X$ is closed, let
\begin{align*}
I_{\mathcal A}(E)&=\{a\in{\mathcal A}:a(x)=0\text{ for }x\in E \}, \\
I_{\mathcal A}^0(E)&=\{a\in{\mathcal A}:a(x)=0\text{ in a nbhd of }E\} \\
\text{ and } J_{\mathcal A}(E)&=\overline{I_{\mathcal A}^0(E)}.
\end{align*}
One says that $E$ is a set of {\it spectral synthesis} for ${\mathcal A}$ if 
$I_{\mathcal A}(E)=J_{\mathcal A}(E)$ (this definition is equivalent to the 
one given in the introduction).

The Banach algebra we will mainly deal with is the projective tensor 
product $V(X,Y)=C(X)\hat\otimes C(Y)$, where
$X$ and $Y$ are compact  Hausdorff spaces.
Recall that $V(X,Y)$ (the Varopoulos algebra) 
consists of all functions 
$\Phi\in C(X\times Y)$ which admit a representation 
\begin{equation}\label{equa}
\Phi(x,y)=\sum_{i=1}^{\infty}f_i(x)g_i(y),
\end{equation}
 where $f_i\in C(X)$, $g_i\in
C(Y)$ and $$\sum_{i=1}^{\infty}||f_i||_{C(X)}||g_i||_{C(Y)}<\infty.$$ 
$V(X,Y)$ is a Banach algebra  with the norm
$$||\Phi||_V=\inf\sum_{i=1}^{\infty}||f_i||_{C(X)}||g_i||_{C(Y)},$$
where $\inf$ is taken over all  representations of $\Phi$ in the form
$\sum f_i(x)g_i(y)$ (shortly, $\sum f_i\otimes g_i$) satisfying the 
above conditions (see \cite{var}).
We note that $V(X,Y)$ is a semi-simple regular Banach algebra with spectra
$X\times Y$. 

For $B\in V(X,Y)'$ and $F\in V(X,Y)$, define $FB$ in $V(X,Y)'$
by $\langle FB,\Psi\rangle=\langle B,F\Psi\rangle$. Define the support
of $B$ by
$$\supp (B)=\{(x,y)\in X\times Y\mid FB\ne 0 \text{ whenever } 
F(x,y)\ne 0\}.$$
Then it is known that for a closed set $E\subseteq X\times Y$,
$$J_{V(X,Y)}(E)^{\perp}=\{B\in V(X,Y)'\mid \supp (B)\subseteq E\}$$
and hence 
$E$ is a 
set of spectral synthesis for $V(X,Y)$ if $I_{V(X,Y)}(E)^{\perp}=
\{B\in V(X,Y)'\mid \supp (B)\subseteq E\}$, i.e., if
$$\langle B,F\rangle=0$$
for any $B\in V(X,Y)'$, $\supp(B)\subseteq E$, and any $F\in V(X,Y)$ 
vanishing on $E$.
Any element of $V(X,Y)'$ can be identified with a 
bounded bilinear form $\langle B,f\otimes g\rangle=B(f,g)$ on
$C(X)\times C(Y)$ which we also call  a bimeasure.

We will need also to consider the class of all functions $\Phi$ on 
$X\times Y$ representable in the form (\ref{equa}) (i.e.
$\Phi(X,Y)=\sum_{i=1}^{\infty}f_i(x)g_i(y)$, where $f_i\in C(X)$, 
$g_i\in C(Y)$) with
$$\sup_x\sum|f_i(x)|^2<\infty,\quad 
\sup_y\sum|g_i(x)|^2<\infty$$ (with the pointwise convergence of the series).
It is called the
extended Haagerup tensor product (\cite{er}) of $C(X)$ and $C(Y)$ and we 
will denote it by
$C(X)\hat\otimes_{eh}C(Y)$. Clearly $V(X,Y)\subset C(X)\hat\otimes_{eh} C(Y)$.
The inclusion is strict, moreover $C(X)\hat\otimes_{eh}C(Y)$ contains some
discontinuous  functions. Indeed,
let $f(x)\in C({\mathbb R})$ such that $|f(x)|\leq 1$,
$f(x)=0$ for any $x\in(-\infty,1]\cup[3/2,+\infty)$ and $f(x)=1$
on the interval $[1+\varepsilon, 3/2-\varepsilon]$, $\varepsilon$ being small
enough. Setting $f_k(x)=f(2^kx)$ and $u(x,y)=\sum f_k(x)\bar{f_k}(y)$, 
we obtain $\sup\sum |f_k(x)|^2=1$ and therefore 
$u(x,y)\in C(X)\hat\otimes_{eh}C(Y)$. However, $u(x,x)=\sum |f_k(x)|^2$ does 
not converge to zero as $x\to 0$ while $u(x,0)=u(0,y)=0$,
i.e. $u(x,y)$ is not continuous in $(0,0)$. On the other hand any function in
$C(X)\hat\otimes_{eh} C(Y)$ is separately continuous and hence it is 
continuous at all points apart of a set of first category.

The following theorem connects operator synthesis and synthesis with respect to
the Varopoulos algebra $V(X,Y)$. 
Let $M(X)$, $M(Y)$ be the spaces of finite Borel measures on $X$ and $Y$ 
respectively.
\begin{theorem1}\label{v}
If a closed set $E\subseteq X\times Y$ is a set of synthesis with respect 
to any pair of measures $(\mu,\nu)$, $\mu\in M(X)$, $\nu\in M(Y)$, then $E$ is 
synthetic with respect to $V(X,Y)$. 
\end{theorem1}
\begin{proof} Assume that $E$ is not a set of spectral synthesis for the
algebra $V(X,Y)$. Then there exists a bimeasure $B$, $\supp (B)\subseteq E$ 
and  $F\in V(X,Y)$, $F\chi_E=0$,  such that $\langle B, F\rangle\ne0$. By the 
Grothendieck theorem, there exist measures $\mu\in M(X)$ and $\nu\in M(Y)$ 
and a constant $C$ such that
\begin{equation}\label{gro}
|\langle B,f\otimes g \rangle|=|B(f,g)|\leq  
C||f||_{L_2(X,\mu)}||g||_{L_2(Y,\nu)}.
\end{equation}
Since $V(X,Y)$ can be densely embedded into $L_2(X,\mu)\hat\otimes L_2(Y,\nu)$,
it follows from (\ref{gro}) that the linear functional 
$\Phi\mapsto \langle B,\Phi\rangle$ defined on $V(X,Y)$ can be extended to a 
continuous linear functional on $L_2(X,\mu)\hat\otimes L_2(Y,\nu)$.
Therefore, there exists an operator $T\in B(L_2(X,\mu),L_2(Y,\nu))$ such that
$$\langle B,\Phi\rangle=\langle T,\Phi\rangle,$$
 the left hand side being  the pairing in the sense of duality between
$V(X,Y)$ and $V(X,Y)'$ and the right hand side is the pairing in the sense
of duality between $L_2(X,\mu)\hat\otimes L_2(Y,\nu)$ and 
$B(L_2(X,\mu),L_2(Y,\nu))$.

We have to  prove that $T$ is supported in
$E$.
Since $E$ is closed, for every closed sets $\alpha$, $\beta$ such that
$(\alpha\times\beta)\cap E=\emptyset$, there exist open sets 
$\alpha_0\supset\alpha$, $\beta_0\supset\beta$ such that 
$\overline{\alpha_0}\times\overline{\beta_0}$ does not intersect $E$.
For every functions $f\in C(X)$, $g\in C(Y)$ which are equal to zero outside
the set $\alpha_0$ and $\beta_0$ respectively, we have
$(Tf,g)=\langle T,f\otimes g\rangle=\langle B,f\otimes g\rangle=0$.
Since each function in $L_2(X,\mu)$ ($L_2(Y,\nu)$) which is zero a.e.\ outside
$\alpha$ ($\beta$) can be approximated
by continous functions vanishing outside $\alpha_0$ 
($\beta_0$ respectively),  we obtain $Q_{\beta}TP_{\alpha}=0$. 
By the regularity of measures $\mu$ and $\nu$ it 
follows that this is true for any Borel sets $\alpha$, $\beta$. 
\end{proof}
\begin{cor1}\label{finw}
Suppose that $\varphi_i:X\mapsto Z$ and $\psi_i:Y\mapsto Z$, 
$i=1,\ldots,n$, are continuous functions from  compact metric
spaces $X$ and $Y$ to an ordered compact metric space $Z$.
Then the set 
$E=\{(x,y)\mid \varphi_i(x)\leq \psi_i(y), 
i=1,\ldots,n\}$
is a set
of synthesis with respect to the algebra $V(X,Y)$.
\end{cor1}
\begin{proof} This  follows from Theorems~\ref{fw},\ref{v}.   
\end{proof}
This corollary yields the theorem of Drury on synthesizability of 
``non-triangular'' sets, which are sets of width two (see \cite{D}).

We will see that  the converse of Theorem~\ref{v} is false in general. 
\begin{lemma1}\label{rectangular}
If $E\subseteq X\times Y$ is a set of synthesis with respect to a pair
of finite measures then so is its intersection with any measurable
rectangle.
\end{lemma1}
\begin{proof} Let $\mu\in M(X)$, $\nu\in M(Y)$, let $K\times S$ be a 
measurable rectangle in $X\times Y$, 
let $T\in B(L_2(X,\mu), L_2(Y,\nu))$ and $F\in \Gamma(X,Y)$ be such that
$\text{supp }T\subseteq E\cap (K\times S)\subseteq \text{null }F$.
Then $T=Q_STP_K$ and $\text{supp }T\subseteq E$. Moreover, the function 
$F'(x,y)=\chi_K(x)\chi_S(y)F(x,y)$ belongs to $\Gamma(X,Y)$ and 
vanishes on $E$. Since
$E$ is a set of synthesis, we obtain 
$$\langle T,F\rangle=\langle Q_STP_K,F\rangle=\langle T,F'\rangle=0,$$
finishing the proof.
\end{proof}

\begin{proposition1}\label{failure}
There exist a closed set $E\subseteq X\times Y$ and a pair $(\mu,\nu)$ of finite
measures on $X$ and $Y$ such that $E$ is set of synthesis 
in $V(X,Y)$, but not of $\mu\times\nu$-synthesis.
\end{proposition1} 

\begin{proof}
It will be  sufficient to
find a closed set $E\subseteq X\times Y$ and a closed rectangle $K\times S$ in 
$X\times Y$ such that $E$ is synthetic with respect to $V(X,Y)$ but not
$E\cap (K\times S)$. In fact, if $E$ were a set of synthesis with respect to 
any pair
of finite measures we would obtain, by Lemma~\ref{rectangular}, that so 
would be  its intersection 
with any measurable rectangle and, by Theorem~\ref{v}, the intersection 
$E\cap (K\times S)$ would be synthetic for $V(X,Y)$.
The construction of the set $E$ is a modification of the Varopoulos example described
in the proof of Theorem~\ref{synfw}.

Let $X$, $Y$ be compact metric spaces and let $G\subset X\times Y$ be 
a non-synthetic set with respect to $V(X,Y)$.  
Let $I$ denote the unit interval $[0,1]$ and $d((x,y),G)$ be the  
distance between $(x,y)$ and $G$. In
$(X\times I)\times Y$ consider the set
$$E=\{((x,t),y)\in (X\times I)\times Y\mid d( (x,y),G)\leq t\}.$$
Then $E$ is a set of synthesis with respect to $V(X\times I,Y)$.
To see this take a function  $F((x,t),y)=\sum_{k=1}^{\infty}f_k(x,t)g_k(y)$
in $V(X\times I,Y)$  such that
\begin{equation}\label{ineq1}
\sum_{k=1}^{\infty} \sup|f_k(x,t)|^2\sum_{k=1}^{\infty} 
\sup|g_k(y))|^2< \infty
\end{equation}
 and $\text{null }F\supseteq E$,
and consider $F_n((x,t),y))=F((x,t+1/n),y)$, $n\in {\mathbb N}$.
Clearly, $F_n$ vanishes on 
$$E_n=\{((x,t),y)\in (X\times I)\times Y\mid d( (x,y),G)< t+1/n\},$$
an open set containing the set $E$.
Now $$F_n((x,t),y)-F((x,t),y)=\sum_{k=1}^{\infty}
(f_k(x,t+1/n)-f_k(x,t))g_k(y)$$ and 
$$||F_n((x,t),y)-F((x,t),y)||_V\leq
\sum_{k=1}^{\infty}
\sup|(f_k(x,t+1/n)-f_k(x,t)|^2\sum_{k=1}^{\infty}\sup |g_k(y)|^2.$$
Fix $\varepsilon>0$. By (\ref{ineq1}) one can find 
$K>0$ such that 
$\sum_{k=K+1}^{\infty}
\sup|(f_k(x,t+1/n)-f_k(x,t)|^2<\varepsilon$. Since all $f_k$, $k=1,\ldots,K$,
are continuous on the compact $X\times I$, they are uniformly continuous.  
Therefore there exists $N>0$ such that,  for any $n\geq N$, we have
$\sup|f_k(x,t+1/n)-f_k(x,t))|<\sqrt{\varepsilon/K}$,  
$k=1,\ldots K$. This yields
 $\sum_{k=1}^{K}
\sup|(f_k(x,t+1/n)-f_k(x,t)|^2<\varepsilon$ and 
$$\sum_{k=1}^{\infty}
\sup|(f_k(x,t+1/n)-f_k(x,t)|^2<2\varepsilon,$$ showing $F_n\to F$ as 
$n\to\infty$ in $V(X\times I,Y)$.

Consider now $$E^*=E\cap ((X\times\{0\})\times Y)=
\{(x,0),y)\in (X\times I)\times Y\mid (x,y)\in G\}.$$
Our goal is to show that $E^*$ is not synthetic in $V(X\times I,Y)$. 
Given a function $\Phi(x,y)=\sum_{k=1}^{\infty} f_k(x)g_k(y)\in V(X,Y)$, 
$\text{null }\Phi\supseteq G$, consider
$F((x,t),y)=\Phi(x,y)$ in $V(X\times I,Y)$. Assume that $E^*$ is synthetic. 
Then $F$ can be approximated in $V(X\times I,Y)$ by functions $F_n((x,t),y)$
which vanish on neighbourhoods of $E^*$. This implies that $\Phi$ can be 
approximated by $F_n((x,0),y)$ in $V(X,Y)$. Clearly, each $F_n((x,0),y)$ 
vanishes on a nbhd of $G$.
By arbitrariness of $\Phi$, we obtain that $G$ is a set of synthesis, 
contradicting our assumption. 
\end{proof}

Thus the sets of universal (independent on the choice of measures) 
operator synthesis form a more narrow class than the sets of spectral 
synthesis. It is of interest to clarify which known classes it includes.

\vspace{0.1cm}

A closed set $E\subseteq X\times Y$ is called ``a set without true bimeasure''
(SWTB, for brevity) if  any bimeasure concentrated on 
$E$ is a measure. It is clear that any such set is a set 
of spectral synthesis in $V(X,Y)$. 
\begin{proposition1}
A closed set without true bimeasures  is a set of universal operator
synthesis.
\end{proposition1}
\begin{proof}
Let $\mu\in M(X)$, $\nu\in M(Y)$ and let $E$ be a closed set without
true bimeasure.  Consider $T\in B(L_2(X,\mu),L_2(Y,\nu))$
such that $T$ is supported in $E$. It defines a bimeasure $B_T$ by
$(Tu,\bar{v})=B_T(u,v)$, where $u\in C(X)$ and $v\in C(Y)$. 
Moreover, $\supp (B_T)\subseteq E$. By the condition of the theorem, there 
exists a measure $m\in M(X\times Y)$ such that $\supp (m)\subseteq E$ and
\begin{equation}\label{measure}
(Tu,\bar{v})=\int u(x)v(y)dm(x,y), 
\end{equation} 
for every $u\in C(X)$, 
$v\in C(Y)$. 

Let $F(x,y)=\sum_{n=1}^{\infty} u_n(x)v_n(y)\in C(X)\hat\otimes_{eh} C(Y)$ and
let $F_k(x,y)=\sum_{n=1}^{k} u_n(x)v_n(y)$, 
$E_k(x,y)=\sum_{n=k+1}^{\infty} |u_n(x)|^2+|v_n(y)|^2$.
Then $E_k(x,y)\to 0$, $k\to\infty$, for every $(x,y)\in X\times Y$,
$$|F(x,y)-F_k(x,y)|\leq E_k(x,y)$$
and therefore $F_k(x,y)\to F(x,y)$, $k\to\infty$, everywhere on $X\times Y$.
Moreover, $|F_k(x,y)|\leq E_0(x,y)$ and $E_0(x,y)$ is integrable 
over $m$, as $m$ is finite.
Thus, by  the theorem on majorized convergence, 
$$\int F_k(x,y)dm(x,y)\to \int F(x,y)dm(x,y).$$
On the other hand, $||F-F_k||_\Gamma\leq\int E_k(x,y)d\mu(x)d\nu(y)$ and
 $\int E_k(x,y)d\mu(x) d\nu(y)\to 0$, which imply
$||F-F_k||_{\Gamma}\to 0$  and $\langle T,F_k\rangle\to\langle T,F\rangle$ 
as $k\to\infty$. 
 
We now obtain the equality
$$\langle T,F\rangle =\int F(x,y)dm(x,y),
\quad  F\in C(X)\hat\otimes_{eh} C(Y).$$ 
Since $m$ is supported  in $E$, this gives  $\langle T, F\rangle=0$ with  $F$ vanishing on $E$.

Consider now $F\in \Gamma(X,Y)$, $\text{null }F\supseteq E$. 
Then there exist $f_i\in L_2(X,\mu)$,
$g_i\in L_2(Y,\nu)$ such that $F(x,y)=\sum_{i=1}^{\infty} f_i(x)g_i(y)$
(m.a.e.) and $\sum_{i=1}^{\infty}||f_i||^2_{L_2}
\sum_{i=1}^{\infty}||g_i||^2_{L_2}<\infty$. Given $\varepsilon>0$, we can 
find compact sets $X_{\varepsilon}\subseteq X$, $Y_{\varepsilon}\subseteq Y$ 
such that $\mu(X\setminus X_{\varepsilon})<\varepsilon$, $\nu(Y\setminus Y_{\varepsilon})
<\varepsilon$ and
\begin{eqnarray*}
\sum_{i=1}^{\infty}|f_i(x)|^2<C_{\varepsilon}, \quad x\in 
X_{\varepsilon},\quad
\sum_{i=1}^{\infty}|g_i(y)|^2<C_{\varepsilon}, \quad y\in  
Y_{\varepsilon},
\end{eqnarray*}
Moreover, we can assume that $f_i$, $g_i$ are continuous by the Lusin theorem
so that the restriction $F_{\varepsilon}$ of $F$ to 
$X_{\varepsilon}\times Y_{\varepsilon}$ belongs to
$C(X_{\varepsilon})\hat\otimes_{eh} C(Y_{\varepsilon})$. Clearly, if
$E$ is a set without true bimeasure, so is $E\cap (X_{\varepsilon}\times
Y_{\varepsilon})$. If now $T\in B(L_2(X,\mu), L_2(Y,\nu))$ is supported in
$E$ then $\text{supp }Q_{Y_{\varepsilon}}TP_{X_{\varepsilon}}\subseteq
E\cap (X_{\varepsilon}\times
Y_{\varepsilon})$ and 
$$\langle Q_{Y_{\varepsilon}}TP_{X_{\varepsilon}}, F_{\varepsilon}\rangle=0.$$
Letting $\varepsilon\to 0$, we obtain $\langle T,F\rangle=0$.
\end{proof}

We can say even more about sets without true bimeasures:  they are
operator solvable (see Definition~2.2).
In the following lemma $(X,\mu)$, $(Y,\nu)$ are finite measure spaces as in 
section~\ref{section2}
\begin{lemma1}\label{solv}
Let $E\subseteq X\times Y$ be a pseudo-closed set. If for any
$\e>0$, there exist $X_{\e}\subseteq X$, $Y_{\e}\subseteq Y$,
$\mu(X\setminus X_{\e})<\e$, $\nu(Y\setminus Y_{\e})<\e$ such that
$E\cap (X_{\e}\times Y_{\e})$ is synthetic in $X_{\e}\times Y_{\e}$ then 
$E$ is synthetic in $X\times Y$.
\end{lemma1}
\begin{proof}
Let $\Phi\in\Gamma(X,Y)$ vanish on $E$. Fix $\e>0$ and set $\Phi_{\e}(x,y)=
\Phi(x,y)\chi_{X_{\e}}(x)\chi_{Y_{\e}}(y)$. 
Clearly, $\Phi_{\e}$ vanishes on $E_{\e}=E\cap (X_{\e}\times Y_{\e})$. Since
$E_{\e}$ is a set of synthesis in $X_{\e}\times Y_{\e}$
there exists $\tilde\Phi_{\e}\in \Gamma (X_{\e},Y_{\e})$ vanishing
in a neighbourhood of $E_{\e}$ such that 
$$||\Phi_{\e}-\tilde\Phi_{\e}||_{\Gamma}<\e.$$
Extending $\tilde\Phi_{\e}$ by zero to the whole space $X\times Y$ we get
a function vanishing on a neighbourhood of $E$ and
\begin{eqnarray*}
||\Phi(x,y)-\tilde\Phi_{\e}(x,y)||_{\Gamma}=
||\Phi(x,y)-\Phi_{\e}(x,y)+\Phi_{\e}(x,y)
-\tilde\Phi_{\e}(x,y)||_{\Gamma}\leq\\
||\Phi(x,y)(\chi_{X_{\e}}(x)\chi_{Y_{\e}}(y)-1)||_{\Gamma}+||\Phi_{\e}(x,y)
-\tilde\Phi_{\e}(x,y)||_{\Gamma}<\\
||\Phi||_{\Gamma}||(\chi_{X_{\e}}\chi_{Y_{\e}}-1)||_{\Gamma}+\e\leq 
||\Phi||_{\Gamma}(\e\nu(Y)+\e\mu(X))+\e,
\end{eqnarray*}
giving the statement.
\end{proof}

\begin{proposition1}
If a closed set $E\subseteq X\times Y$ has a property that any its closed subset
is a set of operator synthesis with respect to a pair $(\mu,\nu)$ of
regular finite measures then $E$ is operator solvable.  In particular,
any set without true bimeasure is operator solvable.
\end{proposition1}
\begin{proof}
Using the regularity of measures, one can easily show that
for any pseudo-closed subset $K\subseteq E$ and any
 $\e>0$ there exists a Borel
rectangle
$X_{\e}\times Y_{\e}$ with  $\mu(X\setminus X_{\e})<\e$, 
$\nu(Y\setminus Y_{\e})<\e$ such that $K\cap (X_{\e}\times Y_{\e})$ is closed.
The statement now follows from Lemma~\ref{solv}. 
\end{proof}

\begin{remark}\rm
In \cite{var}, Varopoulos  
established a deep connection between the algebra 
$V(G)=C(G)\hat\otimes C(G)$ and the Fourier algebra $A(G)$
of compact Abelian groups $G$. Using the relationships he showed that
a closed set $E\subseteq G$ is a set of spectral synthesis for $A(G)$ if
and only if the diagonal set $E^*=\{(x,y)\in G\times G\mid x+y\in E\}$
is a set of spectral synthesis for $V(G)$. 
Recently the same result was proved for non-Abelian compact groups in
\cite{sptu} using the established there connection between $A(G)$ and
the Haagerup tensor product $C(G)\hat\otimes_h C(G)$ which is the Varopoulos 
algebra, renormed. 
An analogous result for sets of operator synthesis in $G\times G$ was obtained
in \cite{F} for locally compact Abelian groups $G$ and in \cite{sptu} for
compact non-Abelian groups $G$. Namely, 
a closed set $E\subseteq G$ is a set of spectral synthesis for $A(G)$
if and only if $E^*$ is a set of operator synthesis with respect to the 
Haar measure (for the reverse statement, synthesizability with respect 
to all pairs 
of finite measures is not required, as in Theorem~\ref{v}).
Using a method  similar to one in Proposition~\ref{failure} one can construct
a set of synthesis $E$ and a pair finite measures $(\mu,\nu)$ such that  
$E^*$ is not $\mu\times\nu$-synthetic.
\end{remark}

\section{Operator-Ditkin sets and union of synthetic sets}

In  the classical harmonic analysis one studies special so-called Ditkin  
(or Wiener-Ditkin or Calderon) sets. If ${\mathcal A}$ is a unital semisimple 
regular commutative Banach algebra with spectrum $X$ then a closed set 
$E\subseteq X$ is called Ditkin set if $u\in \overline{uI_{\mathcal A}^0(E)}$
for every $u\in I_{\mathcal A}(E)$ (see the beginning of the previous section
for the notations).  
An analogue of such sets can be introduced for the space $\Gamma(X,Y)=
L_2(X,\mu)\hat\otimes L_2(Y,\nu)$.
Here we will make use of a space similar to $C(X)\hat\otimes_{eh}C(Y)$
Let
$$
V^{\infty}(X,Y)=L^{\infty}(X,\mu)\otimes^{w^*h}L^{\infty}(Y,\nu)$$
where $\otimes^{w^*h}$ denotes the weak* Haagerup tensor product of 
\cite{BS}.  
$V^{\infty}(X,Y)$ can be identified with a space of functions 
$w:X\times Y\to{\mathbb C}$ which admit 
a representation
$w(x,y)=\sum_{i=1}^\infty\varphi_i(x)\psi_i(y)$,
where $\varphi_i\in L^{\infty}(X,\mu)$, $\psi_i\in L^{\infty}(Y,\nu)$ and
such that the series
$\sum_{i=1}^\infty|\varphi_i|^2$ and $\sum_{i=1}^\infty|\psi_i|^2$ converges
almost everywhere to functions in $L^{\infty}(X,\mu)$ and $L^{\infty}(Y,\nu)$.
As elements in $V^{\infty}(X,Y)$ these functions are defined
up to a marginally null set.

We say that a complex valued function $w$ on $X\times Y$ is a
multiplier of $\Gamma(X,Y)$ if for any $\omega\in\Gamma(X,Y)$,
$(s,t)\mapsto w(s,t)\omega(s,t)$ defines an element  of 
$\Gamma(X,Y)$. One can show that $w$ defines 
a bounded linear operator $m_w$ on $\Gamma(X,Y)$
and two multipliers $w$ and $w'$ satisfy $m_w=m_{w'}$ if
$w=w'$ marginally almost everywhere.  We say that $w$ and $w'$ are
equivalent if $m_w=m_{w'}$.  
It was proved in \cite{sptu} (and in other terms in \cite{peller,smith}) 
that the space of multipliers of
$\Gamma(X,Y)$ coincides with $V^{\infty}(X,Y)$. If measures $\mu$, $\nu$
are finite, we also have $V^{\infty}(X,Y)\subset\Gamma(X,Y)$.
For a pseudo-closed set $E$ denote  $\Psi_{00}(E)=\{F\in V^{\infty}(X,Y):F=0 
\text{ on a neighbourhood of } E\}$.

\begin{definition1}
We say that a pseudo-closed set $E\subseteq X\times Y$ is $\mu\times\nu$-Ditkin 
if
$f\in \overline{f\Psi_{00}(E)}$ for any $f\in\Phi(E)$, i.e.
if for any $f\in\Phi(E)$ there exists a sequence $\{g_n\}\in\Psi_{00}(E)$ such
 that
$$||g_n\cdot f-f||_{\Gamma}\to 0\text{ as }n\to\infty.$$
\end{definition1}

Clearly, every $\mu\times\nu$-Ditkin set is $\mu\times\nu$-synthetic.

We will now study a  question how $\mu\times\nu$-Ditkin and 
$\mu\times\nu$-synthetic sets behave under forming unions.
If $G$ is a locally compact abelian group it is known that
the union of two Ditkin sets in $ X(A(G))$ (the space of characters of the 
Fourier algebra $A(G)$) is Ditkin. Whether the union of two spectral
sets in $A(G)$ is spectral is one of the unsolved problems in harmonic 
analysis. If we knew that any spectral set is a Ditkin set the question would
be answered affirmatively since a union of two Ditkin sets is again a 
Ditkin set (see \cite{benedetto} for survey of this).
Another known result about unions is that if $E$, $F$ are closed subsets 
in $X(A(G))$ such that their
intersection is a Ditkin set then their union is spectral if and only if
so are the sets $E$, $F$ (see \cite{warner}). The result was also generalised
to $A(G)$, where $G$ is an arbitrary locally compact group. 
We will prove a similar statement for $\mu\times\nu$-Ditkin and
$\mu\times\nu$-synthetic sets. In what follows we write simply
Ditkin and synthetic sets, if no confusion arise.

If $f\in\Gamma(X,Y)$ denote by
$\supp (f)=cl_w\{(x,y)\in X\times Y:f(x,y)\ne 0\}$, where $cl_w$ indicates
the $pseudo-closure$. 

\begin{theorem1}
The union of two $\mu\times\nu$-Ditkin sets is a $\mu\times\nu$-Ditkin set.
The union of $\mu\times\nu$-Ditkin set and  a $\mu\times\nu$-synthetic set
is $\mu\times\nu$-synthetic.
\end{theorem1}
\begin{proof}
Suppose $E_1$ and $E_2$ are Ditkin sets, $E=E_1\cup E_2$, $\e>0$,
$f\in \Gamma(X,Y)$ vanishing on $E$. By definition of Ditkin sets there
exist functions  $g_i\in \Psi_{00}(E_i)$ such that
$||f-fg_1||_{\Gamma}<\e/2$ and $||fg_1-fg_1g_2||<\e/2$. If $g=g_1g_2$ then
$g\in V^{\infty}(X,Y)$, $g$ vanishes on a neighbourhood of $E$ and
$||f-fg||_{\Gamma}<\e$.

Let now  $E_1$ be a Ditkin set and $E_2$  synthetic. Then, given $\e>0$ and 
$f\in\Phi(E_1\cup E_2)$, there exist
$g_1\in\Psi_{00}(E_1)$ and $g_2\in\Gamma(X,Y)$ vanishing on a nbhd of 
$E_2$ such that
$||f-fg_1||_{\Gamma}<\e/2$ and $||f-g_2||<e/2||g_1||_{V^{\infty}}$, where
$||g_1||_{V^{\infty}}$ is the norm of the bounded operator on $\Gamma(X,Y)$ 
corresponding to $g_1$. We have that 
$g_1g_2\in \Gamma(X,Y)$ vanishes on a nbhd of $E_1\cup E_2$ and
$||f-g_1g_2||<\e$. 
\end{proof}

\begin{lemma1}\label{dit}
Let $E_1$ and $E_2$ be pseudo-closed subsets of $X\times Y$ whose intersection
is a Ditkin set and let $E=E_1\cup E_2$. Then
$$\Phi_0(E)=\Phi_0(E_1)\cap\Phi_0(E_2).$$
\end{lemma1}
\begin{proof}
Clearly,
$$\Phi_0(E)\subseteq\Phi_0(E_1)\cap\Phi_0(E_2).$$ 
Therefore, we have to prove the reverse inclusion.  We work modulo marginally 
null sets.
Let $f\in\Phi_0(E_1)\cap\Phi_0(E_2)$
 and let  $\delta>0$.
Since $K=E_1\cap E_2$ is a Ditkin set, there is $v\in\Psi_{00}(K)$ such 
that $||vf-f||_{\Gamma}<\delta$.
If $E^c=\cup_{i=1}^{\infty}\alpha_i\times\beta_i$ then
$vf\chi_{\alpha_i\times\beta_i}\in\Phi_0(E)$.
If $(\supp (vf))^c=\cup_{i=1}^{\infty}\gamma_i\times\delta_i$, we have
$vf\chi_{\gamma_i\times\delta_i}=0\in\Phi_0(E)$.

Consider now $E\cap\supp (vf)$ and set $$C_1=E_1\cap\supp (vf),\quad
C_2=E_2\cap\supp (vf).$$ Then $E\cap\supp (vf)=C_1\cup C_2$ and 
$C_1\cap E_2=\emptyset$. 
Hence $C_1\subseteq E_2^c$. Moreover, $E_2^c$ is pseudo-open. 

Let 
$E_2^c=\cup_{i=1}^{\infty}\alpha_i^1\times\beta_i^1$.
We have  $w_i=\chi_{\alpha_i^1\times\beta_i^1}\in \Phi_0(E_2)$ vanishes on a
nbhd of $E_2$
and $vfw_i\in \Phi_0(E_1\cup E_2)=\Phi_0(E)$

Similarly, 
$$C_2\subseteq E_1^c=
\cup_{i=1}^{\infty}\alpha_i^2\times\beta_i^2$$ 
and $vfu_i\in\Phi_0(E)$, where 
$u_i=\chi_{\alpha_i^2\times\beta_i^2}$.
We have therefore
$$X\times Y=E^c\cup(\supp (vf))^c\cup C_1\cup C_2\subseteq\cup_{i=1}^{\infty}
\tilde\alpha_i\times\tilde\beta_i$$
and $vf\chi_{\tilde\alpha_i\times\tilde\beta_i}\in\Phi_0(E)$.
One  can find $A_{\e}\subseteq X$ and $B_{\e}\subseteq Y$, 
$\mu(X\setminus A_{\e})<\e$ and
$\nu(Y\setminus B_{\e})<\e$, such that
$(X\times Y)\cap (A_{\e}\times B_{\e})$ is the union of a finite number
of 
$\{\tilde\alpha_i\times\tilde\beta_i\}$, say first $n$.
Set $v_i=\chi_{\tilde\alpha_i\times\tilde\beta_i}$
and let $h_1=v_1$, $h_2=v_2-h_1v_2,\ldots, h_k=v_k-v_k(h_1+\ldots+h_{k-1})$.
Then $\sum_{i=1}^n h_i=1$ on $A_{\e}\times B_{\e}$, $vfh_i\in\Phi_0(E)$ and
$$vf\chi_{A_{\e}}\chi_{B_{\e}}=vf\sum_{i=1}^n h_i\chi_{A_{\e}}\chi_{B_{\e}}=
\sum_{i=1}^n vfh_i\chi_{A_{\e}}\chi_{B_{\e}}\in\Phi_0(E).$$
Taking now $\e\to 0$ we get $vf\in\Phi_0(E)$ and
$f\in\Phi_0(E)$. 
\end{proof}

\begin{theorem1}
Let $E_1$ and $E_2$ be pseudo-closed subsets of $X\times Y$ whose intersection
is a Ditkin set, and let $E=E_1\cup E_2$. Then $E$ is $\mu\times\nu$-synthetic 
if and
only if both $E_1$ and $E_2$ are $\mu\times\nu$-synthetic.
\end{theorem1}

\begin{proof} Assume first that $E_1$ 
and $E_2$ are synthetic.
We have $$\Phi(E)
=\Phi(E_1)\cap\Phi(E_2)=\Phi_0(E_1)\cap\Phi_0(E_2)=\Phi_0(E).$$
The last equality is due to Lemma~\ref{dit}

To prove the reverse statement we note first that
$\Phi(E)=\Phi_0(E)=\Phi_0(E_1)\cap\Phi_0(E_2)$. On the other hand
$\Phi(E)=\Phi(E_1)\cap\Phi(E_2)$ and we get 
$$\Phi_0(E_1)\cap\Phi_0(E_2)=\Phi(E_1)\cap\Phi(E_2)$$
Take now $f\in\Phi(E_1)$.
Since 
$f\in\Phi(E_1\cap E_2)$ and $K=E_1\cap E_2$ is a Ditkin set, given $\delta>0$
there exists
$v\in\Psi_{00}(K)$ such that $||vf-f||<\delta$. Arguing  as in the proof of
Lemma~\ref{dit} we have
$$E_1^c\cup(\supp (vf))^c=\cup_{i=1}^{\infty}\alpha_i\times\beta_i$$
so that $vf\chi_{\alpha_i\times\beta_i}\in\Phi_0(E_1)$.
Let
$F=E_1\cap\supp (vf)$. We have $F\cap E_2=\emptyset$ and $F\subseteq E_2^c$. 
Then for any $\e>0$ we can find
$X_{\e}\subseteq X$, $Y_{\e}\subseteq Y$, $\mu(X\setminus X_{\e})<\e$ and
$\nu(Y\setminus Y_{\e})<\e$ such that
$$F\cap (X_{\e}\times Y_{\e})\subseteq\cup_{i=1}^{n}\gamma_i\times\delta_i\subseteq
E_2^c.$$
We can choose the rectangles $\gamma_i\times\delta_i$ to be disjoint.

Set $w=\sum_{i=1}^n\chi_{\gamma_i\times\delta_i}$. We have
$vf-vfw$ vanishes on $\cup_{i=1}^{n}\gamma_i\times\delta_i$.
Then
$$X_{\e}\times Y_{\e}\subseteq (F\cup E_1^c\cup(\supp (vf))^c) \cap 
(X_{\e}\times Y_{\e})\subseteq (\cup_{i=1}^{n}\gamma_i\times\delta_i)\cup
(\cup_{i=1}^{\infty}\alpha_i\times\beta_i)$$
and $vf(1-w)\chi_{\cup_{i=1}^{n}\gamma_i\times\delta_i}=0\in\Phi_0(E_1)$,
$vf(1-w)\chi_{\alpha_i\times\beta_i}\in\Phi_0(E_1)$.
As before we can conclude that $$vf(1-w)\chi_{X_{\e}\times
Y_{\e}}\in\Phi_0(E_1).$$
But $vfw\in\Phi_0(E_2)\subseteq\Phi(E_2)$, so $vfw\in\Phi(E_2)\cap
\Phi(E_1)=\Phi_0(E_2)\cap\Phi_0(E_1)$ and therefore $vfw\in\Phi_0(E_1)$. 
Since $(vf-vfw)\chi_{X_{\e}\times Y_{\e}}$ belongs to
$\Phi_0(E_1)$ we get $vf\chi_{X_{\e}\times Y_{\e}}\in\Phi_0(E_1)$. Since
$\e$ and $\delta$ are arbitrary, $vf\in\Phi_0(E_1)$ and $f\in \Phi_0(E_1)$, i.e.
$\Phi_0(E_1)=\Phi(E_1)$. Similarly,
$\Phi_0(E_2)=\Phi(E_2)$.

\end{proof}

\end{document}